\documentclass[11pt]{article}
\usepackage{flafter,amsmath,amssymb,latexsym,psfrag,graphicx,color,indentfirst}
\usepackage[margin=2.5cm]{geometry}

\usepackage[colorlinks,
linktocpage=true,
linkcolor=black,
citecolor=black,
urlcolor=black,
anchorcolor=black
]{hyperref}

\usepackage[numbers,sort&compress]{natbib}
\usepackage{appendix}
\usepackage{multirow}
\usepackage{caption}
\usepackage{subfigure}
\usepackage{graphicx}
\usepackage[normalem]{ulem}

\DeclareGraphicsExtensions{.eps,.pdf,.jpg,.png}
\allowdisplaybreaks

\newtheorem{theorem}{Theorem}[section]

\newtheorem{lemma}[theorem]{Lemma}
\newtheorem{remark}[theorem]{Remark}

\numberwithin{equation}{section}
\numberwithin{figure}{section}

\begin{document}
\setlength\arraycolsep{2pt}
\date{\today}
\title{Recovering both the wave speed and the source function in a time-domain wave equation by injecting contrasting droplets}
\author{Soumen Senapati$^1$, Mourad Sini$^1$, Haibing Wang$^{2,3}$
\\$^1$RICAM, Austrian Academy of Sciences, A-4040, Linz, Austria \qquad\qquad\quad\quad
\\E-mail: soumen.senapati@ricam.oeaw.ac.at and mourad.sini@oeaw.ac.at
\\$^2$School of Mathematics, Southeast University, Nanjing 210096, P.R. China\;
\\$^3$Nanjing Center for Applied Mathematics, Nanjing 211135, P.R. China\qquad\,
\\E-mail: hbwang@seu.edu.cn
}

\maketitle
\begin{abstract}
Dealing with the inverse source problem for the scalar wave equation, we have shown recently that we can reconstruct the space-time dependent source function from the measurement of the wave, collected at a single point $x$ for a large enough interval of time, generated by a small scaled droplets, enjoying large contrasts of its bulk modulus, injected inside the domain to image. Here, we extend this result to reconstruct not only the source function but also the variable wave speed. Indeed, from the measured waves, we first localize the internal values of the travel-time function by looking at the behavior of this collected wave in terms of time. Then from the Eikonal equation, we recover the wave speed. Second, we recover the internal values of the wave generated only by the background (in the absence of the small droplets) from the same measured data by inverting a Volterra integral operator of the second kind. From this reconstructed wave, we recover the source function at the  expense of a numerical differentiation.  
\bigskip

\end{abstract}

\section{Introduction and statement of the results}\label{sec_intro}
\subsection{Introduction}
Let $\Omega$ be a bounded and smooth domain in $\mathbb R^3$ and $T>0$ is some positive real number. Denote by $\rho_0(x)$ and $k_0(x)$ the mass density and bulk modulus of the medium, respectively. Assume that $\rho_0(x)\equiv \rho_0$ is a positive constant in $\mathbb R^3$, but $k_0(x)$ is constant in $\mathbb R^3\setminus\overline\Omega$ and variable in $\Omega$ with positive values. Then the acoustic wave speed in the medium is given by 
\begin{equation}\label{speed_c0}
c_0(x):=\sqrt{\frac{k_0(x)}{\rho_0}}.
\end{equation}
%
%Also the mass density of the medium is denoted by $\rho_0(x)$ which is a positive function in $\mathbb R^3$ admitting a constant value say $\rho_0$ in $\mathbb R^3 \setminus \overline\Omega$. 

Let $J(x,t)$ be a source function which is compactly supported in $\Omega \times (0,T)$, and $v(x,\,t)$ be the corresponding total field generated by the inhomogeneous medium $\Omega$ when probed by the source function $J(x,t)$ (in absence of any initial source). That is, $v$ satisfies the initial value problem
\begin{equation}\label{model_v}
\begin{cases}
c_0^{-2}(x) v_{tt} - \Delta v = J(x,t)  & \textrm{in}\;\mathbb R^3\times (0,\,T), \\
v|_{t=0}=0,\; v_t|_{t=0}=0  &  \mathrm{in}\;\mathbb R^3.
\end{cases}
\end{equation}
The inverse problem that one can immediately formulate, amounts to recover the variable wave speed $c_0$ as well as the source function $J$ from some boundary observation of the wave-field $v$. We will be explicit about our measurement in due time. But before that, let us briefly review the existing literature related to the inverse problem mentioned above. For known wave speed, the determination of source term from a single lateral boundary data has been widely studied since the pioneering work \cite{Bukhgeim_Klibanov_1981} by Bukhge\u{\i}m and Klibanov which mainly deals with the uniqueness part of the  problem and employs Carleman estimate in a crucial way. This method has been subsequently modified in the important works  \cite{Immanuvilov_Yamamoto_2001_CPDE,Yamamoto_1999_JMPA} to address the stability issue of the same problem. All these works are formulated in a bounded domain set-up and assume the wave speed to be constant. Moreover, the source function there is incompletely separated into its spatial and temporal components i.e. $J(x,t) = f(x) R(x,t)$ with $f$ being the only unknown part of $J$ and $R(x,0)$ should admit a strict sign condition throughout $\Omega$. We also mention the work \cite{Bellassoued_2004} where the sound speed is stably recovered in the absence of internal source function $J$, from a single measurement in terms of Dirichlet and Neumann data. Apart from the full Neumann measurement case, there are numerous other works in this direction which treats the same inverse problem but with different boundary data. For instance, the works \cite{Immanuvilov_Yamamoto_2001_IP} and \cite{Bellassoued_Yamamoto_2006_JMPA} consider measurements coming from some suitable interior values of $v$ and from the Neumann observation of $v$ made on arbitrary part of the boundary respectively. We refer to the monograph \cite{Bellassoued_Yamamoto_book} and references therein for an extensive overview on this approach. 
\bigskip

There is a slightly different but quite rich class of inverse problems in connection to ours where one attempts to simultaneously determine the initial source $f(x)$ and sound speed $c_0(x)$ from the Dirichlet data of the wave-field. Said differently, one asks to determine the pair $(c_0,f)$ from the boundary measurement of $v$ which satisfies the IVP \eqref{model_v} but with the changes $J(x,t) \equiv 0$ and $v(x,0) = f(x)$. This problem comes up in photoacoustic tomography (PAT) and thermoacoustic tomography (TAT) which have significant applications in clinical biology and medicine, see for instance \cite{Kruger_Kiser_Reinecke_Kruger,Wang}. Under the assumption that the wave speed is known to be a constant or smooth and non-trapping, the works \cite{Agranovsky_Kuchment_Kunyansky_2009,Finch_Rakesh_2009,Stefanov_Uhlmann_2009,Stefanov_Uhlmann_2013_Transactions} study the recovery of $f$, whereas the work \cite{Stefanov_Uhlmann_2013_Transactions} discusses the recovery of the wave speed $c_0$ for known $f$ under suitable geometric assumption on $\Omega$. The linearized problem of recovering $c_0$ and $f$ simultaneously is quite unstable as shown in \cite{Stefanov_Uhlmann_2013_IPI}. Nonetheless, there are a few works which address the problem of recovering both $c_0$ and $f$ under some conditions.  For admissible class of $(c_0,f)$, the quantity $c_0^{-2}f$ is uniquely determined in \cite{Liu_Uhlmann_2015} where the unknown part of $c_0^{-2}f$ is assumed to be harmonic or free from one spatial variable, leading to the recovery of constant $c_0$ and positive $f$. This argument was further carried on in \cite{Knox_Moradifam_2020}. Very recently, simultaneous recovery of $c_0$ and $f$ was treated in \cite{Kian_Uhlmann_2023} under a monotonicity assumption of $c_0$.  
\bigskip

In our previous work \cite{S-W2022}, we proposed a different approach to solve this inverse problem. The idea is to inject a small droplet (or a bubble) $D$ into $\Omega$ which enjoys a high contrast of its bulk modulus (and eventually its mass density for the bubble, see below for more discussion about these scales) and then measure the generated wave at a point outside or on the boundary of $\Omega$. To be precise, let $D=z+a B\subset\Omega$, where $a$ is a small parameter, $B\subset\mathbb R^3$ is a bounded and smooth domain containing the origin and $z\in\mathbb R^3$.  Define
\begin{equation}\label{def_eps}
k(x):=
\begin{cases}
k_0(x) & \mathrm{in}\; \mathbb R^3\setminus D,\\
k_1 & \mathrm{in}\; D,
\end{cases}
\end{equation}
where $k_1$ is a positive constant with 
\begin{equation}\label{k_1_assumption}
k_1 \sim a^2 \mbox{ as } a\ll 1.
\end{equation}

Set
$$c(x)=\sqrt{\frac{k(x)}{\rho_0 }}, \quad  c_1=\sqrt{\frac{k_1}{\rho_0}}.
$$
Given the same internal source function $J(x,t)$ as before, we denote by $u(x,\,t)$ the wave-field from the medium $\Omega$ with the injected droplet $D$. Then the total field
$u(x,\,t)$
satisfies the following model:
\begin{equation}\label{model_u}
\begin{cases}
c^{-2}(x) u_{tt} - \Delta u = J(x,t) &  \mathrm{in}\;\mathbb R^3\times(0,\,T),\\
u|_{t=0}=0,\;u_t|_{t=0}=0  &  \mathrm{in}\;\mathbb R^3.
\end{cases}
\end{equation}
This can be seen as the linearized model, see \cite{C-M-P-T1986}, of its nonlinear analogue derived in \cite{C-M-P-T1985} for the bubbly media.

In \cite{S-W2022}, we have seen that, if $c_0$ is a known constant, we can reconstruct the source term $J$ from the data given by \begin{equation} u(x, t):=u(x; z, t),\end{equation} collected on a single point $x$ located on $\partial \Omega$ and on a large enough band of time, then  move the location $z$ of the droplet in $\Omega$.  The proposed method is reconstructive. In addition, in terms of dimensionality, the inverse problem is not overdetermined as we use 4D data, 1D in time and 3D for the injected droplets in $\Omega$, to recover a 4D function (1D in time and 3D in space). This approach has its roots back to the works \cite{C-C-S2018, D-G-S2021, A-S2021, A-S2022}. In particular, in \cite{D-G-S2021}, we have shown, in the time-harmonic regime, how we can recover the wave speed, actually the both the mass density and the bulk coefficient, using the resonant character of the injected droplets. In addition, in \cite{A-S2022}, we showed how we can recover both the acoustic and the optic properties of the medium from the Photo acoustic field, generated by injected plasmonic droplets, measured on a single point outside the domain to image $\Omega$. 
\bigskip

In this work, we extend this approach to reconstruct both the speed of propagation $c_0$ and the source term $J$. 
For this, we use as data:
\begin{equation}
  u(x, t):=u(x; z, t) \mbox{ and } v(x, t),
\end{equation}
collected on a single point $x$ located on $\partial \Omega$ and on a large enough band of time. We repeat the measurements of $u(x, t):=u(x; z, t)$ by  moving the location $z$ of the droplet in $\Omega$.

%\begin{remark}\label{First-remark} 
The mathematical modeling described above is related to the ultrasound imaging modality using contrast agents given either by bubbles filled in with gas or droplets filled in with liquids. As these contrast agents are used for imaging in liquids, it is natural to assume that the mass density of the liquid inside the droplets is not varying much from the one of the surrounding liquid. However, the bulk modulus is taken to be very small as compared to the one of the surrounding. Under certain scales, of the form (\ref{k_1_assumption}), these droplets generate local resonances related to eigenvalues of the Newtonian operator stated on its domain.  
Regarding the bubbles, the situation is different. As they are filled in with gas, then it is natural to assume that both their mass density and the bulk modulus are small as compared to the ones of the surrounding liquid.  Under certain scales, these bubbles generate one local resonance, namely the Minnaert resonance frequency which is related to the isolated eigenvalue $\frac{1}{2}$ of the double layer operator stated on the boundary of the bubble. The appearance of these local resonances is key in the analysis. Ultrasound imaging modalities using bubbles or droplets  as contrast agents is well known see \cite{C-F-Q2009, F-M-S2003, I-I-F2018, Q-C-F2009, Q2007, Z-L-L-S-W19}. The generation of the bubbles inside the region of interests can be done at least with two ways. The first way is using cavitation by imposing high pressure, see \cite{Sheeran-Dayton}. The second way is by injecting droplets into that region and then heating them so that the liquid inside evaporates and generate gas inside them, see \cite{Sheeran-Dayton}. In the modeling we propose, we assume that the injection of the droplets is done in isolation. Such injection in isolation of the contrast agents are reported in laboratory experiments, see for instance \cite{H-Z-A2023, Z-L-L-S-W19}.

%\item If in the model (\ref{speed_c0})-(\ref{model_v})-(\ref{def_eps}), we replace $\rho$ by the permittivity $\epsilon$ and $k$ by $\mu^{-1}$, where $\mu$ is the permeability, we get the electromagnetism model in the TE ( or TM) regimes. The condition (\ref{k_1_assumption}) translates into a high contrast of the permeability (or eventually the permittivity) which models dielectric nanoparticles. Therefore, the results presented in this work can be translated to the electromagnetism in the 2D TE or TM regimes (modulo an adaptation of the actual approximation formulas to the 2D-setting).
%\end{enumerate}
%\end{remark}

\subsection{Statement of the results}
Let us set
\begin{equation*}
w(x,\,t):=u(x,\,t)-v(x,\,t),
\end{equation*} 
which is the difference of the wave-fields, before and after injecting the droplet $D$. We notice 
\begin{align}\label{diff eqn}
\begin{cases}
   c_0^{-2}(x) w_{tt} - \Delta w = -\frac{q(x)}{c_0^2(x)} u_{tt} & \textrm{in}\;\mathbb R^3\times (0,\,T), \\
w|_{t=0}=0,\; w_t|_{t=0}=0  &  \mathrm{in}\;\mathbb R^3,
\end{cases}
\end{align}
where $ q(x) = \frac{c_0^2(x)}{c^2(x)} -1$ in $\mathbb R^3$ which is supported in $D$ and of order $a^{-2}$ (in point-wise sense). 

\vspace{2mm}

Let us also point out that the non-constant wave speed $c_0(\cdot)$ gives rise to a Riemannian metric with arc length $d\zeta$ given by
\begin{equation}\label{metric}
d\zeta = \left( \sum_{i=1}^3 c_0^{-2}(x) (dx_i)^2 \right)^{1/2},
\end{equation}
which is related to the travel-time function introduced in Section \ref{Structure of the Green function}. In all of our analysis, we make an assumption on the regularity of the geodesic lines induced by the metric \eqref{metric}, which ensures 
\begin{align}\label{geodesic assumption}
   & \textit{Any pair of points } x,\,y \textit{ in } \Omega \textit{ can be connected with a single geodesic line } \Gamma(x,\,y) \textit{ of the metric } \nonumber\\ 
   & \eqref{metric} \textit{ belonging to } \Omega \textit{ which extends to the boundary of $\Omega$.}
\end{align}
We assume that the source function $J$ is supported in $\Omega \times (0,T)$ and belongs to the function space $H_{0,\sigma}^{13}\left(0,T; L^2(\Omega)\right)$ which is defined in Section \ref{Preliminaries}. 

\begin{theorem}\label{outside domain_4}
   Assuming that the source function $J$ admits the above regularity, $c_0 \in C^{15}(\overline\Omega)$ and the condition \eqref{geodesic assumption} is fulfilled, we have the following asymptotic expansion 
    \begin{align}\label{Main-approximation}
        w(x,t) & = - \sum_{n=1}^{\infty} \frac{\sigma(x,z)c_1}{4 \pi \lambda^{3/2}_n\, c_0(z) |\zeta(x,\,z)|} \left(\int_D e_n(y) \, dy\right)^2 \int_{0}^{t-\zeta(x,z)} \sin\left[\frac{c_1}{\sqrt{\lambda_n}}(t-\tau-\zeta(x,z))\right] \, v(z,\tau)\, d\tau \nonumber\\
    & \quad - \sum_{n=1}^{\infty} \frac{c_1}{\lambda_n^{3/2}} \left(\int_D e_n(y) \, dy\right)^2 \int_{0}^t  \int_0^{t -\tau}\sin\left(\frac{c_1}{\sqrt{\lambda_n}} (t -\tau-s) \right) g(x,s;z) \,v(z,\tau) \, ds\, d\tau \nonumber \\  
    & \quad - v(z,t-\zeta(x,z)) \, \sum_{n=1}^{\infty} \frac{\sigma(x,z)}{4 \pi \lambda_n\, c_0(z) |\zeta(x,\,z)|} \left(\int_D e_n(y) \, dy\right)^2 \nonumber \\
    & \quad - \sum_{n=1}^{\infty} \frac{1}{\lambda_n}  \left(\int_D e_n(y) \, dy\right)^2  \int_{0}^t g(x,t-\tau;z) v(z,\tau)\, d\tau \, + O(a^{2}) ,
    \end{align}
    which holds point-wise in space and time in $\left(\overline{\Omega} \setminus V\right) \times (0,T)$, with $D\subset\subset V \subset\subset \Omega$. 
 The functions $\zeta(\cdot,\cdot),\ \sigma(\cdot,\cdot)$ and $g(\cdot,\cdot;\cdot)$  will be introduced in Section \ref{Preliminaries}, whereas $\{\lambda_n,e_n\}_{n\in\mathbb{N}}$ denotes the eigen-system of the Newtonian operator on $L^2(D)$ which is defined in \eqref{Newtonian_defn}. 
\end{theorem}
\begin{remark}  As it shown in Section \ref{Preliminaries}, we have
\begin{align*}
    c_1 \sim a, \quad \lambda_n \le C a^2, \quad \textrm{and } \quad \left\vert \int_D e_n(x)\, dx \right\vert \le C \, a^{3/2}, \quad n\in\mathbb{N},
\end{align*}
for some $C>0$, independent of $n\in\mathbb{N}$ and $a$.
Therefore, we notice that the first four terms in the expansion of Theorem \ref{outside domain_4} are of order $a$ in point-wise sense w.r.t space and time variables. In light of this, these four terms serve as the dominant terms in the expansion of Theorem \ref{outside domain_4}. A special consideration should be given to the eigenfunctions with non-zero averages, since only the averages of the eigenfunctions appear in Theorem \ref{outside domain_4}.  For spherically shaped $D$, we have nice asymptotic properties and structure theorems for these eigenfamily which can be also used to check the validity of the infinite series in \eqref{Main-approximation}. This aspect has been discussed in detail in Section \ref{recovery of speed and source}. That being said, we now show that the infinite series \eqref{Main-approximation} makes sense even for arbitrary shaped $D$. In this case, we use an observation from \cite{Challa_Mantile_Sini_2020} which says that $u\in H^1(D)$ iff $\sum\limits_{n=1}^\infty \frac{\langle u, e_n\rangle^2}{\lambda_n} < \infty$, where $\langle \cdot,\cdot\rangle$ denotes the standard $L^2$ inner-product. Taking $u = 1$ in $D$, we have
$$\sum\limits_{n=1}^\infty\frac{1}{\lambda_n} \left(\int_D e_n(y)\, dy\right)^2 < \infty,$$  
which ensures the convergence of third and fourth terms in \eqref{Main-approximation}. From the integration by parts,
\begin{align*}
    \int_{0}^{t-\zeta(x,z)} \sin\left[\frac{c_1}{\sqrt{\lambda_n}}(t-\tau-\zeta(x,z))\right] \, v(z,\tau)\, d\tau & = \, \frac{\sqrt{\lambda_n}}{c_1} \,v(z, t-\zeta(x,z)) \\
     - \frac{\sqrt{\lambda_n}}{c_1} &\, \int_{0}^{t-\zeta(x,z)}\cos\left[\frac{c_1}{\sqrt{\lambda_n}}(t-\tau-\zeta(x,z))\right] \, v_t(z,\tau)\, d\tau
\end{align*}
and therefore, we have
\begin{align*}
    & \left\vert \sum_{n=1}^{\infty} \frac{\sigma(x,z)c_1}{4 \pi \lambda^{3/2}_n\, c_0(z) |\zeta(x,\,z)|} \left(\int_D e_n(y) \, dy\right)^2 \int_{0}^{t-\zeta(x,z)} \sin\left[\frac{c_1}{\sqrt{\lambda_n}}(t-\tau-\zeta(x,z))\right] \, v(z,\tau)\, d\tau \nonumber\right\vert \\
    & \preceq \sigma(x,z) \sup_{0\le \tau\le t}|v(z,\tau)|\, \sum\limits_{n=1}^\infty\frac{1}{\lambda_n} \left( \int_D e_n(y)\, dy\right)^2 < \infty.
\end{align*}
As a consequence, the first infinite series in \eqref{Main-approximation} is convergent. Similarly, we can establish the validity of second infinite series in \eqref{Main-approximation}. 
\bigskip

In \cite{S-W2022}, the background was taken to be homogeneous, therefore the terms involving $g$ in (\ref{Main-approximation}) disappear. In addition, in \cite{S-W2022}, we have used the approximation for a large enough time $t$ so that, due to the support, in space and time, of the source function $J(\cdot, \cdot)$ the forth term in (\ref{Main-approximation}) disappears as well. In conclusion, (\ref{Main-approximation}) reduces to the one derived in \cite{S-W2022} if the background is homogeneous and $t$ is large enough. 
\end{remark}

From the measurement of $w(x,\cdot)$, we now proceed to determine the sound speed and source function. Our strategy here is to recover the internal travel-time function and the initial wave-field $v(\cdot,\cdot)$. Having achieved this, the knowledge of $\zeta(x,\cdot)$ then determines the wave speed $c_0(\cdot)$ via the Eikonal equation \eqref{Eikonal}. Whereas we explicitly use $c_0(\cdot)$ and $v(\cdot,\cdot)$ to recover the source function from \eqref{model_v}. This is based on the following result.

\begin{theorem}\label{speed_source}
     Let us assume that $x\in \partial{\Omega}, \, c_0 \in C^{15}(\overline\Omega),\ J\in H_{0,\sigma}^{13}\left(0,T; L^2(\Omega)\right)$ is compactly supported in $\Omega \times (0,T)$ and the condition  \eqref{geodesic assumption} holds. Then there exists a linear and invertible map $ \mathbb{A}: L^2(0,T) \rightarrow L^2(0,T)$ such that
    \begin{align}\label{expansion_v}
        v(z,t) = \mathbb{A} w(x,\cdot)(t+\zeta(x,z)) + O(a) ,
    \end{align}
which holds point-wise in space variables and $L^2$-sense (but also point-wise) in time variables. 
Under the assumption that $J(y,\cdot)$ is not identically zero for a.e. $y\in\Omega$, there exists a dense set $\mathcal{W}$ in $\Omega$, such that for $z\in \mathcal{W}$, $\mathbb{A}^{-1} v(z,t)= 0$, when $t\le \zeta(x,z)$ and $ \mathbb{A}^{-1} v(z,t) \neq 0$, when $t > \zeta(x,z)$. 
\end{theorem}
%
%\subsection{A discussion about the results} 
The operator $\mathbb{A}$ is defined in Section \ref{recovery of speed and source}.
 There, we see that $\mathbb{A}w(x,\cdot)(t) \simeq O(1)$, to be understood point-wise in time and thus it is the dominant term in the right hand side of \eqref{expansion_v}. It shows that the wave-field $v(z,\cdot)$ can be recovered from $w(x,\cdot)$ in $L^2$-sense modulo an error term $O(a)$. The time-regularity of $v$ help us obtain $v$ in the whole of $(0,T)$. Working only with the dominant part of the measurement, the travel-time function $\zeta(x,z)$ can be uniquely identified as the time-level where the graph of the function $\mathbb{A}^{-1} v(z,\cdot - \zeta(x,z))$ experiences a jump for the first time \footnote{Due to the time-regularity of $v$ and definition of $\mathbb{A}$ in Section \ref{recovery of speed and source}, we can give sense to the jump of the function $\mathbb{A}^{-1} v(z,\cdot - \zeta(x,z))$ w.r.t the time variable.}. Therefore, this time level can be located (upto an error) from the observation of $w(x, t)$ with respect to time $t\in (0, T)$ (with a single point $x$ on $\partial\Omega$). Once the function $\zeta(x,\cdot)$ is known, we use the Eikonal equation \eqref{Eikonal} to determine the wave speed $c_0$. Having said that, we need to exploit the regularity of $\zeta(x,\cdot)$ to recover $c_0$ since we can determine the function $\zeta(x,\cdot)$  at most in  a dense set. A detailed discussion on this issue  is presented in Section \ref{remarks}.

\vspace{2mm}

 In addition, from (\ref{expansion_v}), we recover $v(z, t), \; z \in \Omega$ and $t \in (0, T)$. Therefore, knowing $c_0$ and $v(\cdot, \cdot)$, we can recover the source function $J(\cdot, \cdot)$ from (\ref{model_v}) via a numerical differentiation.   We note that the representation of measured wave-field and the reconstruction scheme of initial wave-field involve an infinite series. The later is due to the expansion \eqref{inverse_series expansion}. From the viewpoint of real applications, it is desirable to obtain the stability and accuracy results of the inversion scheme. Hence, one should truncate the series expansion of these wave-fields according to the noise level present in the data. The numerical consideration of our reconstruction scheme will be addressed in a forthcoming work.

\vspace{2mm}

Compared to the known literature, as cited earlier, we provide a reconstruction method that uses the above described data (a $4$ dimensional manifold) to recover both the wave speed and the source term (also a $4$ dimensional manifold).  This makes it a competitive result. An attractive feature of our reconstruction method lies in the notion of `locality', which is desired in many practical applications. It indicates that, to recover the medium properties of a region, one needs to place the droplets in a small neighbourhood of that specific region only. Hence, it is economical from the viewpoint of practitioners. Let us also mention that the source term $J(\cdot, \cdot)$ could be replaced by an initial term, as for photoacoustics and thermoacoustics. A similar analysis would justify the reconstruction of both the wave speed and the initial source term. Regarding the photoacoustic tomography problem, the reader can find in \cite{A-S2022} a justification of a reconstruction scheme that uses plasmonic nanoparticles as contrast agents. With such data, a simultaneous reconstruction of the wave speed, the mass density and the electric permittivity of the medium was shown. This provides, at once, a full solution of the photoacoustic tomography problem as we recover, at once, both the acoustic and optical properties of the medium to image.  
\bigskip

In terms of the inverse problem, at the analysis level, the difference between the current work, as presented in Theorem \ref{speed_source}, and the one in \cite{S-W2022} can be understood as follows:
\begin{enumerate}
\item In \cite{S-W2022}, we dealt with the reconstruction of the source term $J$ assuming that the background is known and homogeneous (i.e. the wave speed is a known constant). There, we worked with the time measurement at some point $x\in\mathbb{R}^3 \setminus\Omega$ of the wave $u$ generated by the source term in the presence of the small droplet.  Therefore, we did not need to measure the primary wave $v(x, t)$ at the point $x$ (i.e. we measure only after having injected the small droplet $u(x, t)$). The corresponding inversion scheme was based on a constructed  Riesz basis. The construction of this basis was possible for spherically shaped small droplets for which we have nice asymptotic properties of the eigenvalues of the corresponding Newtonian operator.

\item In the current work, in addition to the source function $J$, we also handle an unknown background (i.e. unknown wave speed). For this, we consider early time of measurement. Therefore, we need to measure both $u(x, t)$ and $v(x, t)$ (i.e. $w(x, t):=u(x, t)-v(x, t)$) for a given point $x \in \partial{\Omega}$. 
The definition and the inversion of the operator $\alpha I+  \mathcal{K}$ is possible for spherically shaped small droplets, as it is based on nice properties of the spectrum of the Newtonian operator as well.  However, as stated before, with such measurements we can reconstruct internal values of the travel-time function as well as the internal values of the primary field $v(\cdot, \cdot)$. With these two quantities, we can reconstruct both the wave speed $c_0(\cdot)$ and the source function $J(\cdot, \cdot)$.
\end{enumerate}

 Let us mention that, in the remaining parts of this work, we follow the notation, $f(x) = O(a^r)$ to imply that there is $C>0$ independent of $a$ so that $|f(x)| \le C a^r$ holds true in the domain of definition of $f$. Also, we use the notation $ a\preceq b$ (or, $a \succeq b$) to imply $a \le C \, b$ (or, $ C\, a \ge b$) for some $C>0$.

\section{Preliminaries}\label{Preliminaries}

Let us introduce the appropriate function spaces which will be used throughout this article. For $r\in \mathbb{R}$ and $T \in (0,\infty]$, we define 
\begin{align*}
    H^r_0(0,T) = \{f\rvert_{(0,T)}: \ f \in H^r(\mathbb{R}) \textnormal{ and } f\rvert_{(-\infty,0)} \equiv 0\}.
\end{align*}
        Likewise, one can introduce similar notion of a generalized space consisting of $E$-valued function where $E$ is a Hilbert space and denote it by $H^r_0(0,T; E)$. For $\sigma>0$ and $r\in \mathbb{Z}_{+}$, we define 
\begin{align*}
    H^r_{0,\sigma}(0,T; E) = \left\{ f \in H^r_0(0,T; E); \ \sum_{k=0}^{r} \int_0^T e^{-2\sigma t} \|\partial^k_t f(\cdot,t)\|_{E}^2 \, dt < \infty \right\}.
\end{align*}
It is not difficult to see that $H^r_0(0,T; E) = H^r_{0,\sigma}(0,T; E)$ for $T < \infty$. However, there is  strict inclusion when $T=\infty$ i.e. $H^r_{0,\sigma}(0,T; E) \subsetneq H^r_0(0,T; E)$. In any case, we will follow the notation 
\begin{align*}
    \|f\|^{2}_{H^r_{0,\sigma}(0,T; E)} := \sum_{k=0}^{r} \int_0^T e^{-2\sigma t} \|\partial^k_t f(\cdot,t)\|_{E}^2 \, dt.
\end{align*}
As we have seen in Theorem \ref{outside domain_4}, the asymptotic expansion of the wave-field $w$  crucially relies on the eigen-system of the Newtonian operator $\mathcal{N}_{D}$ defined by
\begin{align}\label{Newtonian_defn}
  \mathcal{N}_{D}(f)(x) = \int_D \frac{f(y)}{4\pi\,|x-y|} \, dy , \quad  \textnormal{ for } x \in D.
\end{align}
We can similarly define the Newtonian operator $\mathcal{N}_{B}$ on $L^2(B).$  It is straight forward to see that if $\tilde e_n$ is a normalized eigenfucntion of $ \mathcal{N}_{B}$ corresponding to the eigenvalue $\tilde\lambda_n$, then 
\[
e_n(x) = \frac{1}{a^{3/2}} \tilde e_n\left(\frac{x-z}{a}\right) , \quad \textnormal{ for } x \in D
\]
is a normalized eigenfunction of $ \mathcal{N}_{D}$ corresponding to the eigenvalue $\lambda_n = a^2 \tilde\lambda_n$. Furthermore, 
\begin{align*}
    \int_D e_n(x) \ dx = a^{3/2} \int_B \tilde{e}_n(x) \ dx, \quad \forall n \in \mathbb{N}.
\end{align*}
 Since, $\{\tilde\lambda_n\}_{n\in\mathbb{N}}$ denotes the sequence of eigenvalues corresponding to $\mathcal{N}_{B}$, we have $\lim\limits_{n\to\infty}\tilde\lambda_n=0$. Therefore, the relation $\lambda_n \le Ca^2, $ holds true for any $n\in\mathbb{N}$ where $C>0$ is independent of $n$ and $a$. From the Cauchy-Schwarz inequality, it is clear that 
\begin{align*}
    \left|\int_B \tilde{e}_n(x) \ dx\right\vert \le |B|^{1/2} \|\tilde{e}_n\|_{L^2(B)} = |B|^{1/2}, \quad n\in\mathbb{N}
\end{align*}
implying $\int_D e_n(x) \ dx \le C\, a^{3/2}$, for $n\in\mathbb{N}$, where $C>0$ does not depend on $n$ and $a$.

\subsection{A-priori estimates}

To study the asymptotic behavior of the wave-field $u(x,\,t)$ near $D\times(0,\,T)$, which solves the IVP \eqref{model_u}, we need to discuss the existence and a-priori estimate results for $u$ and $v$. It follows from a standard argument. However, we sketch the proof in the following for the sake of completeness. We only pursue the case for $u$ and keep track of the small parameter $a$ while deriving its a-priori estimates. The IVP for $v$ is free from the parameter $a$ and can be handled in a similar manner.
\begin{lemma}\label{a-priori estimate}
    For $J\in H^{p}_{0,\sigma}\left(\mathbb{R}_+; L^2(\mathbb{R}^3)\right)$, we have $ u \in H^{p+1}_{0,\sigma}\left((0,T); L^2(\Omega)\right)$ which solves \eqref{model_u} and satisfies
    \begin{align*}
        \|u\|_{H^{p+1-r}_{0,\sigma}\left((0,T); H^r(\Omega)\right)} \le  \frac{C}{a^{r^2-r}} \|J\|_{H^{p}_{0,\sigma}\left(\mathbb{R}_+; L^2(\mathbb{R}^3)\right)}, \quad r \in \{0,1,2\},
    \end{align*}
    where $C>0$ is some constant independent of the parameter $a>0$. 
\end{lemma}
\textbf{Proof.}
    With the causality assumption imposed on our source function $J$, we aim to show that $u$ too is a causal function lying in appropriate space which solves \eqref{model_u}. Now, let us consider the elliptic problem 
\begin{align}\label{Laplace trans eqn}
  - \Delta \tilde{u}(x,s) + \frac{s^2}{c^2(x)} \tilde{u}(x,s) = \hat{J}(x,s), \quad  \textnormal{ for } \Re(s) = \sigma >0,
\end{align}
which can be thought of the Laplace transformed version of  \eqref{model_u}.
Here we followed the convention that $\hat{J}(x,\cdot)$ denotes the Laplace transformation of $J(x,\cdot)$. The unique solution to \eqref{Laplace trans eqn} lying in $H^1(\mathbb R^3)$ can be found by a variational argument after introducing the sesquilinear map $\mathbb{B}[\cdot,\cdot]: H^1(\mathbb R^3) \times H^1(\mathbb R^3) \to \mathbb{C}$ and the antilinear map $\mathbb{J}: H^1(\mathbb R^3) \to \mathbb{C}$ as
\begin{align*}
\mathbb{B}[\phi,\psi] & = \int_{\mathbb R^3} \nabla \phi(x) \cdot \overline{\nabla \psi(x)} \, dx + s^2 \int_{\mathbb R^3} \frac{\phi(x)\overline{\psi(x)}}{c^2(x)}\, dx ,  \\
\mathbb{J}(\psi) & = \int_{\mathbb R^3} \hat{J}(x,s) \overline{\psi(x)} \, dx .
\end{align*}
The antilineraity of $\mathbb{A}$ and sesquilinearity of $\mathbb{B}$ follow from the definitions. We observe that $\mathbb{B}$ may not be coercive. Hence we consider the sesqulinear map $\mathbb{B}[\cdot, s\, \cdot]$ which becomes coercive since 
\begin{align*}
    \Re\left(\mathbb{B}[\phi, s\phi]\right) = \Re\left(\bar{s} \int_{\mathbb{R}^3} |\nabla\phi(x)|^2 \, dx \, + s |s|^2 \int_{\mathbb{R}^3} \frac{|\phi(x)|^2}{c^2(x)} \, dx \right) \succeq \min\{\sigma, \sigma^3\} \|\phi\|^2_{H^1(\mathbb{R}^3)}.
\end{align*}
In view of this, the weak formulation of \eqref{Laplace trans eqn} amounts to finding $\tilde{u}(\cdot,s) \in H^1(\mathbb{R}^3)$ so that
\begin{align*}
    \mathbb{B}[\tilde{u}(\cdot,s), s \phi] = \bar{s} \, \mathbb{J}(\phi), \quad \forall \phi\in H^1(\mathbb{R}^3),
\end{align*}
which is guaranteed by Lax-Milgram theorem. In consequence, we have
\begin{align}\label{pre-coercivity}
 \bar{s} \|\nabla \tilde{u}(\cdot,s)\|^2_{L^2(\mathbb R^3)} + s |s|^2 \left\|c^{-1} \tilde{u}(\cdot,s)\right\|^2_{L^2(\mathbb R^3)} = \bar{s} \int_{\mathbb{R}^3} \hat{J}(x,s) \bar{\tilde{u}}(x,s) \, dx .
\end{align}
Considering only the real parts of \eqref{pre-coercivity} with an use of Cauchy-Schwarz inequality, we can infer 
\begin{align*}
    \sigma \left( \|\nabla \tilde{u}(\cdot,s)\|^2_{L^2(\mathbb R^3)} + |s|^2 \left\|c^{-1} \tilde{u}(\cdot,s)\right\|^2_{L^2(\mathbb R^3)} \right) \le |s| \|c\hat{J}(\cdot,s)\|_{L^2(\mathbb R^3)} \|c^{-1} \tilde{u}(\cdot,s)\|_{L^2(\mathbb R^3)}, 
\end{align*}
which further implies
\begin{align}\label{ineq for a priori}
  \|\tilde{u}(\cdot,s)\|_{L^2(\mathbb R^3)} \le \frac{C}{\sigma |s|} \|\hat{J}(\cdot,s)\|_{L^2(\mathbb R^3)}, \  \|\nabla \tilde{u}(\cdot,s)\|_{L^2(\mathbb R^3)} \le \|\hat{J}(\cdot,s)\|_{L^2(\mathbb R^3)}.
\end{align}
Here $C$ is a positive constant depending only on the uniform upper bound of $c(\cdot)$ which is not related to the scaling parameter $a>0$. Moreover, we can invoke interior $H^2$-regularity results (in spatial variables) to improve the regularity of $\tilde{u}(\cdot,s)$.  After using \eqref{ineq for a priori}, we eventually obtain from \eqref{Laplace trans eqn} the interior regularity estimate  
\begin{align}\label{ineq for a priori_2}
  \nonumber \|D^2\tilde{u}(\cdot,s)\|_{L^2(\Omega)} & \le C \left( \left\|\frac{s^2}{c^2(\cdot)}u(\cdot) - \hat{J}(\cdot,s)\right\|_{L^2(\mathbb R^3)}+ \|\tilde{u}(\cdot,s)\|_{L^2(\mathbb R^3)}\right) \\
  & \le \frac{C}{a^2} |s| \|\hat{J}(\cdot,s)\|_{L^2(\mathbb R^3)}.
\end{align}
The constants $C>0$ above are independent of $a$ and depends only on $\Omega$ and $T$.

Let us now define the inverse Laplace transform of $\tilde{u}(x,\cdot)$ as follows
\begin{align}\label{defn_u}
    u(x,t) := \int_{\sigma- i\infty}^{\sigma+i\infty} e^{st} \tilde{u}(x,s) \, ds  = e^{\sigma t} \int_\mathbb{R} e^{it\mu} \tilde{u}(x,\sigma+i\mu) \, d\mu.
\end{align}
Notice that our definition of $u$ does not rely on the choice $\Re(s) = \sigma$. It can be justified by performing complex integration of $e^{st} \tilde{u}(x,s)$ on a rectangular contour in $\mathbb{C}_{+}$ and then letting off the top and bottom part of that contour to infinity.  Due to the estimate in \eqref{ineq for a priori}, the function $u$ defined in \eqref{defn_u} is causal and admits polynomial bounds in time. For a detailed discussion regarding this, we refer the reader to consult \cite[Chapter 3]{Sayas_2016_book}. Meanwhile, we see that
\begin{align*}
    \mathcal{F}_{t \to \mu} \left(e^{-\sigma t} \partial^k_t u(x,t)\right) = s^k \tilde{u}(x,s), \quad s = \sigma + i \mu,
\end{align*}
where $\mathcal{F}_t$ denotes Fourier transform w.r.t time variable. As a result, we deduce with a use of \eqref{ineq for a priori}  
\begin{align*}
    \|u\|^2_{H^{p+1}_{0,\sigma}\left(\mathbb{R}_+; L^2(\mathbb{R}^3)\right)} & = \int_0^{\infty} e^{-2\sigma t} \sum_{k=0}^{p+1} \|\partial^k_t u(\cdot,t)\|^2_{L^2(\mathbb{R}^3)} \,dt \\
  & \preceq \int_{\mathbb{R}_+} \int_\Omega e^{-2\sigma t} \sum_{k=0}^{p+1} |\partial^k_t u(x,t)|^2 \, dx\, dt \\
  & \preceq \int_{\mathbb{R}^3} \int_{\mathbb{R}} \sum_{k=0}^{p+1} \left|\mathcal{F} \left( e^{-\sigma t} \partial^k_t u\right)(x,t)\right|^2 \, dt\, dx \\
  & \preceq \sum_{k=0}^{p +1} \int_{\sigma+i\mathbb{R}} |s|^{2k} |\tilde{u}(\cdot,s)|^2_{L^2(\mathbb R^3)} \, ds \\
  & \preceq \sum_{k=0}^{p} \int_{\sigma+i\mathbb{R}} |s|^{2k} |\hat{J}(\cdot,s)|^2_{L^2(\mathbb R^3)} \, ds \simeq \|J\|^2_{H^{p}_{0,\sigma}\left(\mathbb{R}_+; L^2(\mathbb{R}^3)\right)}.
\end{align*}
Following similar arguments and utilizing \eqref{ineq for a priori} and \eqref{ineq for a priori_2}, we can also establish
\begin{align*}
 \|u\|_{H^{p}_{0,\sigma}\left(\mathbb{R}_+; H^1(\mathbb{R}^3)\right)} \preceq \|J\|_{H^{p}_{0,\sigma}\left(\mathbb{R}_+; L^2(\mathbb{R}^3)\right)}, \quad \|u\|_{H^{p-1}_{0,\sigma}\left(\mathbb{R}_+; H^2(\Omega)\right)} \le \frac{C}{a^2} \|J\|_{H^{p}_{0,\sigma}\left(\mathbb{R}_+; L^2(\mathbb{R}^3)\right)}.
\end{align*}

Now we show that $u$ defined in \eqref{defn_u} solves the IVP \eqref{model_u}. For that purpose, let us choose $\psi \in H^1(\mathbb{R}^3)$. For a.e. $t$, the definition of $u$ and the weak formulation of \eqref{Laplace trans eqn} imply 
\begin{align*}
    & \int_{\mathbb{R}^3} c_0^{-2}(x) \partial^2_t u(x,t) \overline{\psi(x)} \, dx + \int_{\mathbb{R}^3} \nabla u(x,t) \cdot \overline{\nabla \psi(x)} \, dx  \\
    & = \int_{\sigma- i\infty}^{\sigma+i\infty} \int_{\mathbb{R}^3} e^{st} \left( c_0^{-2}(x) s^2 \tilde{u}(x,s) \, \overline{\psi(x)}  +  \nabla \tilde{u}(x,s) \cdot \overline{\nabla \psi(x)}  \right) \, dx \, ds \\
    & = \int_{\sigma- i\infty}^{\sigma+i\infty} \int_{\mathbb{R}^3} e^{st} \hat{J}(x,s) \overline{\psi(x)} \, dx \, ds \\
    & = \int_{\mathbb{R}^3} J(x,t) \overline{\psi(x)} \, dx.
\end{align*}
Hence the proof is complete.

\subsection{Structure of the Green's function}\label{Structure of the Green function}

 Here, we show the local singularity of the Green's function $G(x,\,t;\,y,\,\tau)$ for $x,\,y\in\Omega$, and then analyze the Lippmann-Schwinger equation for $u(x,\,t)$ in order to derive its asymptotic expansion in $D \times (0,T)$.

Introduce the function  $H(x,\,t;\,y,\,\tau)$ defined by
\begin{equation}\label{v_Green_H}
\begin{cases}H_{tt} - c_0^2(x)  \Delta H = \delta(t-\tau)\,\delta(x-y)  &  \mathrm{in}\;\mathbb R^3\times \mathbb R,\\
H|_{t<\tau}=0  &  \mathrm{in}\;\mathbb R^3.
\end{cases}
\end{equation}
Let $\zeta(x,\,y)$ be a solution to the Eikonal equation
\begin{equation}\label{Eikonal}
c_0^2(x)\,|\nabla_x \zeta(x,\,y)|^2 = 1
\end{equation}
with
\begin{equation*}
 \zeta(x,\,y)=O(|x-y|)  \quad \textrm{as}\; x\to y.
\end{equation*}
The physical meaning of $ \zeta(x,\,y)$ is the traveling time of a signal from $y$ to $x$. 
At the point $y\in \Omega$, the point $x$ can be given with the Riemann coordinates $\xi=(\xi_1,\,\xi_2,\,\xi_3)$, i.e., $x=f(\xi,\,y)$. With respect to $\xi$, the function $f(\xi,\,y)$ has the inverse $\xi=\eta(x,\,y)$ with
\begin{equation}\label{v_g}
\eta(x,\,y)=x-y + O(|x-y|^2) \quad \textrm{as}\; x\to y.
\end{equation}
Assume that $c_0(x)\in C^{9}(\overline\Omega),\,0<c_0(x)\leq C_0$. Then we have
\begin{equation}\label{v_H2}
\zeta(x,\,y)=c_0^{-1}(y)\, |\eta(x,\,y)| = c_0^{-1}(y)\,|x-y| + O(|x-y|^2) \quad \textrm{as}\; x\to y,
\end{equation} 
which can be easily observed from (3.16) and (3.17) in \cite{Romanov1987}.

 We denote the Riemann ellipsoid by $$S(x,\,y,\,t):=\{\xi:\,\zeta(\xi,\,y) + \zeta(\xi,\,x) =t \},$$ which is well defined for $t \geq \zeta(x,\,y)$ and contracts to $\Gamma(x,\,y)$ as $t \to \zeta(x,\,y)$. Let $\mathcal G$ be a set of points $(x,\,t,\,y)$ such that $x,\,y\in\Omega,\,t>\zeta(x,\,y)$ and $S(x,\,y,\,t)\subset \Omega$. Set $\mathcal G_0 = \mathcal G\cup \mathcal G^\prime$, where 
\begin{equation}\label{G'}
\mathcal G^\prime := \{(x,\,t,\,y):\,x,\,y\in\Omega,\,0\leq t \leq \zeta(x,\,y)\}.
\end{equation}

According to Theorem 4.1 in \cite{Romanov1987}, the function $H(x,\,t;\,y,\,\tau)$ can be written as
\begin{align}\label{v_H1}
 H(x,\,t;\,y,\,\tau) =  \frac{\sigma(x,\,y)}{2\pi c_0^{3}(y)} \,\delta\left((t-\tau)^2-\zeta^2(x,\,y)\right) + h(x,\,t-\tau;\,y), \quad t\geq \tau, \,x\not= y, 
\end{align}
where 
\begin{equation*}
\sigma(x,\,y) =  | \nabla_x \eta(x,\,y) |^{1/2}\, \exp\left( -\frac{1}{2} \int_{\Gamma(x,y)} \nabla(c_0^2)\cdot \nabla\zeta \, d\zeta\right),
\end{equation*}
and $h(x,\,t-\tau;\, y)$ is continuous in $\mathcal G$, piecewise continuous in $\mathcal G_0$ and supported for points satisfying $ t > \tau + \zeta(x,y)$. Note that 
%$$ | \nabla_x g(y,\,y) |=1, \quad g(y,\,y)=0.$$ 
%It follows that
%\begin{equation*}
%\sigma_{-1}(x,\,y)=c_0^{-3}(y) + K(x,\,y),
%\end{equation*}
%where $K(x,\,y)$ is a smooth function with $K(x,\,y)=O(|x-y|)$ as $x\to y$. In addition, we have
%
\begin{eqnarray*}
\delta\left((t-\tau)^2 - \zeta^2(x,\,y)\right)  &=& \frac{1}{2 |\zeta(x,\,y)|}\left[ \delta(t-\tau-|\zeta(x,\,y)|) +  \delta(t-\tau+|\zeta(x,\,y)|) \right] \\
&=& \frac{1}{2 |\zeta(x,\,y)|} \delta(t-\tau-|\zeta(x,\,y)|), \quad t\geq \tau,\;x\not=y.
\end{eqnarray*}
Hence, we obtain
\begin{eqnarray}\label{v_H3}
H(x,\,t;\,y,\,\tau)  =  \frac{1}{4\pi} \frac{\sigma(x,y)}{c_0^3(y)\,|\zeta(x,\,y)|}\,\delta(t-\tau -\zeta(x,\,y))  + h(x,\,t-\tau;\,y).
\end{eqnarray}
Here, we would like to point out that the regularity assumption on the wave speed i.e. $c_0 \in C^{15}(\overline\Omega)$ is made to entertain the function $h(x,t;y)$ to be Lipschitz w.r.t $y$ variable inside the characteristic conoid $t > \zeta(x,y)$. For details, consult \cite[Theorem 4.1]{Romanov1987} for the case when $s=2$ and $l=11$.

The Green's function $G$ satisfies
\begin{equation}\label{v_Green2}
\begin{cases}G_{tt} - c_0^2(x) \Delta G = c_0^2(x)\,  \delta(t-\tau)\,\delta(x-y)  &  \mathrm{in}\;\mathbb R^3\times \mathbb R,\\
G|_{t=\tau}=0,\;G_t|_{t=\tau}=0  &  \mathrm{in}\;\mathbb R^3.
\end{cases}
\end{equation}
Then it follows that
\begin{eqnarray*}\label{v_H4}
G(x,\,t;\,y,\,\tau) &=&\int_{\mathbb R}\int_{\mathbb R^3} H(x,\,t;\,z,\,s)\,c_0^2(z)\,  \delta(s-\tau)\,\delta(z-y) \,dzds \nonumber\\
& =& c_0^2(y) \, H(x,\,t;\,y,\,\tau) \nonumber\\
&=&  \frac{1}{4\pi} \frac{\sigma(x,y)}{c_0(y)\,|\zeta(x,\,y)|}\,\delta(t-\tau -\zeta(x,\,y)) + g(x,\,t-\tau;\,y)
\end{eqnarray*}
where we took $g(x,t-\tau;y):= c_0^2(y)\, h(x,\,t-\tau;\,y)$ which evidently is Lipschitz in $y$ variable inside the characteristic conoid $t > \zeta(x,y)$. 
We state this representation as a theorem:
\begin{theorem}  With the assumption \eqref{geodesic assumption} fulfilled and $c_0 \in C^{15}(\overline\Omega)$, we have the following representation of $G$ in terms of its singularities at the second order:
\begin{equation}\label{v_H4_} 
G(x,\,t;\,y,\,\tau) =  \frac{1}{4\pi} \frac{\sigma(x,y)}{c_0(y)\,|\zeta(x,\,y)|}\,\delta(t-\tau -\zeta(x,\,y)) + g(x,\,t-\tau;\,y).
\end{equation}
\end{theorem}

\bigskip

\section{Proof of Theorem \ref{outside domain_4}}\label{sec_asy}

Let us briefly discuss what follows in this section. Broadly speaking, our attempt to derive an asymptotic expansion for the wave-field $u$ away (and near) $D \times (0,T)$ can be divided into two steps. In the first step, we consider the volume integral representation of $u$ and convert it to an ODE involving $u$ and it's Newtonian transform which is defined in \eqref{Newtonian_defn}. The computation here brings along several error terms which are mainly handled by using Lemma \ref{a-priori estimate}. In the second step, we improve our previously obtained $L^2$ bound of $u$ (and its' time-derivatives) on $D$ in terms of the scaling parameter $a$, for each time. To accomplish this, we resort to the use of spectral theory for the Newtonian operator which has been already employed in \cite{S-W2022}. The difference is that the kernel of the operator is not the fundamental solution, as in \cite{S-W2022}, but the Green's function $G$. To handle this difference, we heavily employ the singularity decomposition (\ref{v_H4_}).

In view of \eqref{diff eqn}, we express 
\begin{align*}
    w(x,t) = - V_G \left[\frac{q}{c_0^2} u_{tt}\right](x,t),
\end{align*}
which we call the volume integral representation of $w$. Here $V_G$ denotes the volume operator having kernel as $G(x,t;\cdot,\cdot)$, that is
\begin{align*}
V_G \left[f\right](x,t) = \int_{\mathbb R^3}\int_{\mathbb R} G(x,t;y,\tau) f(y,\tau) \, dy\, d\tau.
\end{align*}
Since $ w(x,t) = u(x,t) - v(x,t)$, we are led to the following integral equation in $u$
\begin{align}\label{vol eqn}
  u(x,t) + V_G \left[\frac{q}{c_0^2} u_{tt}\right](x,t) = v(x,t),
\end{align}
which is known as the Lippmann-Schwinger equation. We analyze \eqref{vol eqn} by retaining only the dominant terms and collecting the rest into an error term. On that note, we recall the structure of $G$ given in \eqref{v_H4_} and obtain
\begin{align*}
   & V_G \left[\frac{q}{c_0^2} u_{tt}\right](x,t) \\
    & = \int_{\mathbb{R}^3} \int_{\mathbb R} \frac{1}{4\pi} \frac{\sigma(x,y)}{c_0(y)\,|\zeta(x,\,y)|}\,\delta(t-\tau -\zeta(x,\,y)) \frac{q(y)}{c_0^2(y)} u_{tt}(y,\tau)\, dy \, d\tau \\
    & \quad + \int_{\mathbb{R}^3} \int_{\mathbb R} g(x,t-\tau;y) \frac{q(y)}{c_0^2(y)} u_{tt}(y,\tau)\, dy \, d\tau \\
    & = I_1(x,t) + I_2(x,t).
\end{align*} 
We now evaluate both the terms $I_1$ and $I_2$. In order to do so, we first notice that
\begin{align}\label{error terms_1}
    \frac{\sigma(x,y)q(y)}{c_0^3(y)\zeta(x,y)} = \frac{q(z)}{c_0^2(z) |x-y|} + \frac{O(a^{-1})}{|x-y|} + O(a^{-2}), \quad \textrm{ as } y \to x.
\end{align}
It can be verified by straight forward calculations starting with the consideration
\begin{align}\label{error terms_2}
   \frac{\sigma(x,y) q(y)}{c_0^3(y)\zeta(x,y)} =  \frac{q(z)}{c_0^2(z) |x-y|} + \sum_{i=1}^{4} \textnormal{err}_{(i)}(x,y) 
\end{align}
where
\begin{align*}
    \textnormal{err}_{(1)}(x,y) &:=  \frac{q(y)\sigma(x,y) - q(x)}{c_0^3(y)\zeta(x,y)}, \\
    \textnormal{err}_{(2)}(x,y) &:=  c_0^{-2}(y) q(x)\left(\frac{1}{c_0(y)\zeta(x,y)} - \frac{1}{|x-y|} \right), \\
    \textnormal{err}_{(3)}(x,y) &:=  \frac{q(x)(c_0^{-2}(y) - c_0^{-2}(x))}{|x-y|}, \\
    \textnormal{err}_{(4)}(x,y) &:=  \frac{c_0^{-2}(x)q(x) - c_0^{-2}(z) q(z)}{|x-y|}.
\end{align*}
Here we utilize the relation \eqref{v_H2} and regularity of the functions $q(\cdot)\sigma(\cdot,\cdot), \  q(\cdot) c_0^{-2}(\cdot)$ and $c_0^{-2}(\cdot)$ in $D$. To be precise, the first three error terms in RHS of \eqref{error terms_2} which are $\textnormal{err}_{(i)}(x,y),\, i=1,..,3$, will produce the $O(a^{-2})$ term in \eqref{error terms_1}. Proving this for the third error term requires only the smoothness of $c_0^{-2}(\cdot)$, that is
\begin{align*}
    \left\vert \textnormal{err}_{(3)}(x,y) \right\vert \preceq \frac{O(a^{-2})|x-y|}{|x-y|} = O(a^{-2})
\end{align*}
For $y$ sufficiently close to $x$, we use \eqref{v_H2} to observe that
\begin{align}\label{travel-time est}
   |x-y| \preceq \zeta(x,y). 
\end{align}
With the smoothness of $q(\cdot)\sigma(\cdot, \cdot),\ c_0^{-3}(\cdot)$ and the inequality \eqref{travel-time est} in hand, we argue
\begin{align*}
    \left\vert \textnormal{err}_{(1)}(x,y) \right\vert \preceq \left\vert\frac{q(y)\sigma(x,y) - q(x)}{c_0^3(y)\zeta(x,y)} \right\vert \preceq \frac{O(a^{-2}) |x-y|}{|\zeta(x,y)|} = O(a^{-2}).
\end{align*}
and, in a similar way using \eqref{v_H2}, we obtain
\begin{align*}
    \left\vert \textnormal{err}_{(2)}(x,y) \right\vert & \preceq c_0^{-2}(y)q(x)\left\vert\frac{1}{c_0(y)\zeta(x,y)} - \frac{1}{|x-y|} \right\vert \\
    & \preceq c_0^{-2}(y) q(x)\left\vert\frac{\zeta(x,y)-c_0^{-1}(y)|x-y|}{\zeta(x,y)|x-y|}\right\vert \\
    & \preceq O(a^{-2}) \frac{O\left(|x-y|^2\right)}{|x-y|^2} \\
    & \preceq O(a^{-2}).
\end{align*}
Likewise, the fourth error term in \eqref{error terms_2} can be estimated as 
\begin{align*}
    \left\vert \textnormal{err}_{(4)}(x,y) \right\vert \preceq \left\vert \frac{c_0^{-2}(x)q(x) - c_0^{-2}(z)q(z)}{|x-y|} \right\vert \preceq O(a^{-2}) \frac{|x-z|}{|x-y|} \preceq \frac{O(a^{-1})}{|x-y|}.
\end{align*}
Since the term $q$ is supported in $D$, we write
\begin{align}\label{breaking Lippmann-Schwinger}
   I_1(x,t) & = \frac{1}{4\pi} \int_{D} \frac{\sigma(x,y)q(y)}{c^3_0(y)\,|\zeta(x,\,y)|} u_{tt}(y,t-\zeta(x,y)) \, dy  \nonumber \\
       & = \frac{1}{4\pi} \int_{D} \frac{\sigma(x,y)q(y)}{c^3_0(y)\,|\zeta(x,\,y)|} \left(u_{tt}(y,t) + \int_{0}^{\zeta(x,y)} \partial^3_t u(y,t-\kappa) \, d\kappa \right)\, dy \nonumber\\
       & = \frac{q(z)}{c^2_0(z)} \int_D \frac{u_{tt}(y, t)}{4\pi|x-y|} \, dy \, + \sum_{i=1}^{4} \textnormal{err}_{(i)}(x,t) ,
\end{align}
where 
\begin{align*}
     \textnormal{err}_{(1)}(x,t) & = \frac{q(z)}{c^2_0(z)}  \int_D \frac{1}{4\pi |x-y|}\int_{0}^{\zeta(x,y)} \partial^3_t u(y,t-\kappa) \, d\kappa dy, \\
     \textnormal{err}_{(2)}(x,t) & = O(a^{-1}) \int_D \frac{u_{tt}(y, t)}{4\pi|x-y|} \, dy,  \\
     \textnormal{err}_{(3)}(x,t) & = O(a^{-1}) \int_D \frac{1}{4\pi |x-y|}\int_{0}^{\zeta(x,y)} \partial^3_t u(y,t-\kappa) \, d\kappa dy, \\
     \textnormal{err}_{(4)}(x,t) & = O(a^{-2}) \int_D u_{tt}(y,t-\zeta(x,y)) \, dy. 
\end{align*}
Keeping the first term as it is, we focus on  the remaining terms from \eqref{breaking Lippmann-Schwinger}. We will shortly see these remaining terms enjoy better $L^2$ norms with respect to the parameter $a$, when compared against the first term. On account of the a-priori estimates of $u$ from Lemma \ref{a-priori estimate}, the relation \eqref{travel-time est} and H\"older's inequality, we notice that
\begin{align}\label{error_1}
 \|\textnormal{err}_{(1)}(\cdot,t)\|_{L^2(D)} & = \left\|\frac{q(z)}{c^2_0(z)}  \int_D \frac{1}{4\pi |x-y|}\int_{0}^{\zeta(x,y)} \partial^3_t u(y,t-\kappa) \, d\kappa \, dy\right\|_{L^2(D)} \nonumber \\
& \preceq |q(z)| \left(\int_D \left| \int_D \frac{1}{4\pi |x-y|}\int_{0}^{\zeta(x,y)} \partial^3_t u(y,t-\kappa) \, d\kappa \, dy\right|^2 \, dx \right)^{1/2} \nonumber \\
& \preceq |q(z)| \left(\int_D \int_D \frac{|D|}{|x-y|^2} \left|\zeta(x,y)\right|^2 \underset{s \in [0,T]}{\textnormal{sup}}\left|\partial^3_t u(y,s)\right|^2  \, dy\, dx\right)^{1/2} \nonumber \\
& \preceq |q(z)| \left( |D|^{2} \int_D \underset{s \in [0,T]}{\textnormal{sup}}\left|\partial^3_t u(y,s)\right|^2 \, dy\right)^{1/2} \nonumber \\
& = O(a) ,
\end{align}
and in a similar way
\begin{align*}
   \|\textnormal{err}_{(3)}(\cdot,t)\|_{L^2(D)} = \left\|O(a^{-1}) \int_D \frac{1}{4\pi |x-y|}\int_{0}^{\zeta(x,y)} \partial^3_t u(y,t-\kappa) \, d\kappa \, dy\right\|_{L^2(D)} \preceq O(a^2).
\end{align*}
We also compute
\begin{align*}
   \|\textnormal{err}_{(2)}(\cdot,t)\|_{L^2(D)} & = \left\| O(a^{-1}) \int_D \frac{u_{tt}(y, t)}{4\pi|x-y|} \, dy \right\|_{L^2(D)} \\
    & = O(a^{-1}) \left( \int_D \left\vert \int_D \frac{u_{tt}(y, t)}{|x-y|} \, dy\right\vert^2  \, dx \right)^{\frac{1}{2}}\\
    & \preceq O(a^{-1}) \|\partial_t^2 u(\cdot,t)\|_{L^\infty(D)}  \left( \int_{ D}  |D| \int_{ D} \frac{1}{|x-y|^2} \, dy \, dx\right)^{\frac{1}{2}} \\
    & \preceq O(\sqrt a) \|\partial_t^2 u(\cdot,t)\|_{L^\infty(D)} \, a^2  \left(\int_{B(0,1)} \int_{B(0,1)}\frac{1}{|x-y|^2} \, dy \, dx\right)^{\frac{1}{2}} \\
    & = O(\sqrt a),
\end{align*}
where Sobolev embedding theorem along with Lemma \ref{a-priori estimate} has been employed to obtain
\begin{align*}
    \|\partial_t^2 u(\cdot,t)\|_{L^\infty(D)} \preceq \|u\|_{H^3_{0,\sigma}(0,T; H^2(\Omega))} \preceq a^{-2}, \quad \textrm{ for a.e. } t \in (0,T).
\end{align*}
Again we take help from Lemma \ref{a-priori estimate} to derive
\begin{align}\label{error_2}
\|\textnormal{err}_{(4)}(\cdot,t)\|_{L^2(D)} & = \left\|O(a^{-2}) \int_D u_{tt}(y,t-\zeta(x,y)) \, dy\right\|_{L^2(D)} \nonumber \\
& = O(a^{-2}) \left( \int_D \left|  \int_D u_{tt}(y,t-\zeta(x,y)) \, dy\right|^2 \, dx\right)^{1/2} \nonumber \\
& \preceq O(a^{-2}) \left(|D|^{2} \int_D \underset{s \in [0,T]}{\textnormal{sup}}\left|\partial^2_t u(y,s)\right|^2 \, dy \right)^{1/2} \nonumber \\
& = O(a).
\end{align}

Now, we focus on estimating the term $I_2$. We recall that $q$ is supported in $D$ and $g$ in $\mathcal G^\prime$, see (\ref{G'}). Therefore, we have
\begin{align*}
   I_2(x,t) = \int_{\mathbb{R}^3} \int_{\mathbb R} g(x,t-\tau;y) \frac{q(y)}{c_0^2(y)} u_{tt}(y,\tau)\, dy d\tau = \int_D \int_0^{t-\zeta(x,y)} g(x,t-\tau;y) \frac{q(y)}{c_0^2(y)} u_{tt}(y,\tau)\, dy d\tau.
\end{align*}
Due to the boundedness of $g$ in $\mathcal G^\prime$ and Lemma \ref{a-priori estimate}, we deduce
\begin{align}\label{error_3}
   \|I_2(\cdot,t)\|_{L^2(D)} & \preceq \left(\int_D\left|\int_D \int_0^{t-\zeta(x,y)} g(x,t-\tau;y) \frac{q(y)}{c_0^2(y)} u_{tt}(y,\tau)\, dy d\tau \right|^2 dx\right)^{1/2} \nonumber \\
   & \preceq |q(z)| \left( |D|^{2} \int_D \underset{s \in [0,T]}{\textnormal{sup}}\left|\partial^2_t u(y,s)\right|^2 \, dy\right)^{1/2} \nonumber \\
   & = O(a).
\end{align}
Summarizing the above calculations, the Lippmann-Schwinger equation \eqref{vol eqn} reduces to 
\begin{align}\label{reduction to Newtonian op}
\begin{cases}
  u(x,t) + \frac{q(z)}{c^2_0(z)} \int_D \frac{u_{tt}(y, t)}{4\pi|x-y|} \, dy \, = f(x,t), \quad \textnormal{ for } (x,t) \in D \times (0,T), \\
  u(x,0) = u_t(x,0) = 0,    
\end{cases}
\end{align}
where the source term $f$ is given by
\begin{align}\label{source}
    f(x,t) := v(x,t) - \sum\limits_{i=1}^{4} \textnormal{err}_{(i)}(x,t) - I_2(x,t) = v(x,t) + O(\sqrt a),
\end{align}
to be understood point-wise sense in time and $L^2$-sense in space. This concludes the discussion of first step hinted at the beginning.

The second step concentrates on improving the $L^2$ bounds for $u$ and its time-derivatives upto third order. In order to do so, we study the problem \eqref{reduction to Newtonian op}. We work with general source function say $f$ in \eqref{reduction to Newtonian op} for future use and thereafter prove the result
\begin{align}\label{bound in New op}
   \|u(\cdot,t)\|_{L^2(D)} \preceq \underset{0 \le \tau \le t}{\textnormal{sup}} \left(\|f(\cdot,\tau)\|_{L^2(D)} + \|f'(\cdot,\tau)\|_{L^2(D)} \right), \quad \textrm{ for } t\in (0,T),
\end{align}
where $u$ solves \eqref{reduction to Newtonian op}. To find $u$, we set
\begin{align*}
  u(\cdot,t) = \sum_{n=1}^{\infty} u_n(t) e_n , \quad f(\cdot, t) =  \sum_{n=1}^{\infty} f_n(t) e_n , \quad \textnormal{ in } L^2(D),
\end{align*}
or equivalently,
\begin{align*}
    u_n(t) = \langle u(\cdot,t), e_n\rangle_{L^2(D)},  \quad f_n(t) = \langle f(\cdot,t), e_n\rangle_{L^2(D)},  \quad \forall t \in (0,T), \ n \in \mathbb{N}.
\end{align*}
Here $\{\lambda_n,e_n(\cdot)\}_{n=1}^{\infty}$ denotes an orthonormal eigensystem of the Newtonian operator $\mathcal{N}_{D}$ on $L^2(D)$.
The coefficients $\{u_n(t)\}_{n=1}^{\infty}$ are to be determined from the following ODE which follows due to the linear independence of $\{e_n\}_{n\in \mathbb{N}}$,
\begin{align*}
\begin{cases}
  u_n(t) + \frac{\lambda_n q(z)}{c_0^2(z)} u^{''}_n(t) = f_n(t), \quad n \in \mathbb{N}, \\
  u_n(0) = u'_n(0) = 0 .
\end{cases}
\end{align*}
The unique solution to this ODE can be explicitly written down as 
\begin{align*}
  u_n(t) = c_0(z) (q(z)\lambda_n)^{-1/2} \int_0^t \sin\left[c_0(z)\,(q(z)\lambda_n)^{-1/2}(t-\tau)\right]\, f_n(\tau)\,d\tau ,
\end{align*}
and performing an integration by parts, we assert that
\begin{align}\label{u_n}
    u_n(t) = -\int_0^t \cos\left[c_0(z)\,(q(z)\lambda_n)^{-1/2}(t-\tau)\right] f'_n(\tau)\,d\tau \, + f_n(t).
\end{align}
After an application of Cauchy-Schwarz inequality to \eqref{u_n}, we obtain
\begin{align*}
    |u_n(t)|^2 \preceq \underset{0 \le \tau \le t}{\textnormal{sup}} \left( |f_n(\tau)|^2 + |f'_n(\tau)|^2 \right), \quad \textrm{ for a.e. }  t \in (0,T).
\end{align*}
From the representation of $u$ above, we see 
\[ 
\|u(\cdot,t)\|^2_{L^2(D)} = \sum\limits_{n=1}^{\infty} |u_n(t)|^2,
\]
and similar equality holds for $f(\cdot,t)$ and $f'(\cdot,t)$. Hence, we arrive at the desired estimate \eqref{bound in New op}. 

Applying the estimate \eqref{bound in New op} in our case i.e. when $f$ is given explicitly in \eqref{source}, we have 
\begin{align*}
    \|u(\cdot,t)\|_{L^2(D)} \preceq \|v(\cdot,t)\|_{L^2(D)} + \sqrt{a} \preceq \sqrt{a}, \quad \textnormal{ a.e. } t \in (0,T),
\end{align*}
which is already an improvement of the bound for $\|u(\cdot,t)\|_{L^2(D)}$ when one applies Lemma \ref{a-priori estimate}. However, we want to improve it upto $O(a^{3/2})$ even for time-derivatives of $u$ upto third order. We can repeat our earlier argument but for $\partial_t^k u$ by taking time-derivatives of $k$-th order to \eqref{reduction to Newtonian op}. In that case, we need to ensure $L^2$ integrability of higher time-derivatives of the source term $f$ in \eqref{reduction to Newtonian op}, which is indeed the case here. For instance, if we try to employ the estimate \eqref{bound in New op} for $\partial_t^3 u$, the source term then consists of functions like $\partial_t^3 v$ and $\partial_t^6 u$. As a consequence, we find
\begin{align}\label{first bound}
   \|\partial_t^k u(\cdot,t)\|_{L^2(D)} \preceq \|\partial_t^k v(\cdot,t)\|_{L^2(D)} + \sqrt a \preceq \sqrt a, \quad \textnormal{ for a.e. } t \in (0,T), \, k \in \{1,..,11\}. 
\end{align}
We can utilize \eqref{first bound} to further improve the error estimates \eqref{error_1}-\eqref{error_3}. To be specific, we see
\begin{align*}
    \sum_{i=1}^{4} \|\textnormal{err}_{(i)}\|_{H^1(0,T;L^2(D))} + \|I_2\|_{H^1(0,T;L^2(D))} \preceq a^{3/2},
\end{align*}
and reconsider the problem \eqref{reduction to Newtonian op}. But now, the source term $f$ satisfies 
\begin{align*}
    \underset{0 \le \tau \le t}{\textnormal{sup}} \left(\|f(\cdot,\tau)\|_{L^2(D)} + \|f'(\cdot,\tau)\|_{L^2(D)} \right) \preceq \|f\|_{H^1(0,T;L^2(D))} \preceq a^{3/2}.
\end{align*}
With the help of \eqref{bound in New op}, we have 
\begin{align*}
   \|u(\cdot,t)\|_{L^2(D)} \preceq a^{3/2} ,\quad \textnormal{ a.e. } t\in (0,T).
\end{align*}
For a.e. $t \in (0,T)$ and $k \in \{1,2,3\}$, we proceed as in \eqref{first bound} and finally obtain
\begin{align}\label{second bound}
    \|\partial_t^k u(\cdot,t)\|_{L^2(D)} \preceq \|\partial_t^k v(\cdot,t)\|_{L^2(D)} + a^{3/2} \preceq a^{3/2}.
\end{align}
Therefore, the error estimates \eqref{error_1}-\eqref{error_3} can be improved once more. As a result, the problem \eqref{reduction to Newtonian op} boils down to its revised form which is
\begin{align*}
     u(x,t) + \frac{q(z)}{c^2_0(z)} \int_D \frac{u_{tt}(y, t)}{4\pi|x-y|} \, dy \, = v(x,t) + O(a^{5/2}), \quad \textnormal{ in } D \times (0,T),
\end{align*}
which holds in point-wise sense w.r.t time and $L^2$-sense w.r.t space variables.

\bigskip

Now, we establish the asymptotic expansion of the wave-field $u$ near $D \times (0,T)$. Shortly after, we also derive the same but for $w$ in $  \left(\overline{\Omega} \setminus V \right)\times (0,T)$ where $D\subset\subset V\subset\subset \Omega$. In that way, we get an idea of the structure of the measured wave-field. Taking $f(x,t) := v(z,t)$ for $x \in D$ in \eqref{bound in New op}  and denoting the corresponding solution by $u_z$, we notice that
\begin{align*}
   f_n(x,t) = v(z,t) \int_D e_n(x) \, dx, \quad \textrm{ for } n \in \mathbb{N}, \, t \in (0,T),
\end{align*}
and thus $u_z$ has the following expression in $D \times (0,T)$:
\begin{align*}
u_z(x,t) = \sum_{n=1}^{\infty} c_0(z) (q(z)\lambda_n)^{-1/2} e_n(x) \int_0^t \sin\left[c_0(z)\,(q(z)\lambda_n)^{-1/2}(t-\tau)\right] v(z,\tau)\, d\tau \int_D e_n(y) \, dy.
\end{align*}
A major objective in introducing the wave-field $u_z$ is that it approximates the original field $u$ modulo the error $O(a^{2})$ to be understood point-wise in time and $L^2$-wise in $D$. Also, the dominant part of $u$, i.e. $u_z$, requires the wave-field $v$ only at the point $z$ where the droplet has been injected. We discuss the details in the following. 

Let us define the error 
\begin{align*}
   E(x,t) = u(x,t) - u_z(x,t), \quad \textrm{ for } (x,t) \in D \times (0,T),
\end{align*}
which satisfies 
\begin{align*}
   E(x,t) + \frac{q(z)}{c^2_0(z)} \int_D \frac{E_{tt}(y, t)}{4\pi|x-y|} \, dy \, = v(x,t) - v(z,t) + O(a^{5/2}) = O(a^{2})
\end{align*}
in point-wise sense in time and $L^2$-sense w.r.t space variables. Here, we make use of Sobolev embedding theorem to gurantee $v(\cdot,t) \in C^{0, \frac{1}{2}}(\overline\Omega)$ since $v(\cdot,t) \in H^2(\Omega)$ and therefore, for a.e. $t\in (0,T)$
\begin{align*}
    \|v(\cdot,t) - v(z,t) \|_{L^2(D)} = \left(\int_{D} |v(x,t) - v(z,t)|^2 \, dx\right)^{1/2} \preceq \left( \int_{D} |x-z| \, dx \right)^{1/2} \preceq a^2. 
\end{align*}
Now, we appeal to \eqref{bound in New op} to draw the conclusion 
\[
u(x,t) = u_z(x,t) + O(a^{2}),\ \quad \textnormal{ for } (x,t) \in D \times (0,T), 
\]
which leads to the asymptotic expansion 
\begin{align}\label{expansion in D}
   u(x,t) & = \sum_{n=1}^{\infty} c_0(z) (q(z)\lambda_n)^{-1/2} e_n(x) \int_0^t \sin\left[c_0(z)\,(q(z)\lambda_n)^{-1/2}(t-\tau)\right] v(z,\tau)\, d\tau \int_D e_n(y) \, dy \nonumber \\
   & \quad + O(a^{2}) \nonumber \\
   & = \sum_{n=1}^{\infty} \frac{c_1}{\sqrt{\lambda_n}} e_n(x) \int_0^t \sin\left[\frac{c_1}{\sqrt{\lambda_n}}(t-\tau)\right] v(z,\tau)\, d\tau \int_D e_n(y) \, dy + O(a^{2}) ,
\end{align}
to be understood point-wise sense in time and $L^2$-sense in $D$. Based on our previous discussion, we can further claim that the time-derivatives upto third order of the error term in \eqref{expansion in D} will be also of order $a^{2}$ (point-wise in time and $L^2$-sense in $D$). 

\bigskip

Now, we derive the asymptotic expansion of $w$ outside $\overline{D}$. Let us first fix $x \in \overline{\Omega}\setminus V$. In order to do so, we use the structure of the fundamental solution $G(x,t;y,\tau)$ given in \eqref{v_H4_} and write
\begin{align}\label{exapnsion_w}
  \nonumber w(x,t) & = - \int_{\mathbb R^3}\int_{\mathbb R} G(x,t;y,\tau) \frac{q(y)}{c_0^2(y)} u_{tt}(y,\tau) \, dy \, d\tau \\
   \nonumber & = - \frac{1}{4\pi} \int_{D} \frac{\sigma(x,y)q(y)}{c^3_0(y)\,|\zeta(x,\,y)|} u_{tt}(y,t-\zeta(x,y)) \, dy \, - \int_{D} \int_{0}^{t} g(x,t-\tau;y) \frac{q(y)}{c_0^2(y)} u_{tt}(y,\tau)\, dy \, d\tau \\
  \nonumber & = - \frac{\sigma(x,z)q(z)}{4 \pi c^3_0(z)\,|\zeta(x,\,z)|} \int_{D} u_{tt}(y,t-\zeta(x,z)) \, dy \, - \frac{q(z)}{c_0^2(z)} \int_{D} \int_{0}^{t} g(x,t-\tau;z)  u_{tt}(y,\tau)\, dy \, d\tau \\
   & \quad + \sum_{i=1}^3 {\textnormal{err}}_{(i)}(t),
\end{align} 
where 
\begin{align*}
    & {\textnormal{err}}_{(1)}(t):= - \frac{1}{4\pi} \int_{D} \left(\frac{\sigma(x,y)q(y)}{c^3_0(y)\,|\zeta(x,\,y)|} - \frac{\sigma(x,z)q(z)}{c^3_0(z)\,|\zeta(x,\,z)|}\right)u_{tt}(y,t-\zeta(x,z)) \, dy,   \\
    & {\textnormal{err}}_{(2)}(t):= - \frac{1}{4\pi} \int_{D} \frac{\sigma(x,y)q(y)}{c^3_0(y)\,|\zeta(x,\,y)|} \left(u_{tt}(y,t-\zeta(x,y)) - u_{tt}(y,t-\zeta(x,z))\right) \, dy, \\
    & {\textnormal{err}}_{(3)}(t):= - \int_{D} \int_{0}^{t} \left(\frac{q(y)g(x,t-\tau;y)}{c_0^2(y)} - \frac{q(z)g(x,t-\tau;z)}{c_0^2(z)}\right)  u_{tt}(y,\tau)\, dy \, d\tau.
\end{align*}
We show that the error terms ${\textnormal{err}}_{(i)}(t),\ i=1,2,3$ are of order $a^2$ to be understood pointwise in time. Using the smoothness of $\zeta(x,\cdot),\, \sigma(x,\cdot)$ and $c_0(\cdot)$, we find
\begin{align}
 \label{reln_outside domain_1}  \frac{\sigma(x,y)q(y)}{c^3_0(y)\,|\zeta(x,\,y)|} & = \frac{\sigma(x,z)q(z)}{c^3_0(z)\,|\zeta(x,\,z)|} + O(a^{-2})|y-z|, \quad \textnormal{ for } y \in D, \\
  \label{reln_outside domain_3} u_{tt}(y,t-\zeta(x,y)) & = u_{tt}(y,t-\zeta(x,z))  + \int_{\zeta(x,z)}^{\zeta(x,y)} \partial^3_t u(y,t-\kappa) \, d\kappa , \quad \textnormal{ for } y \in D,
\end{align}
and also we have
\begin{align*}
  {\textnormal{err}}_{(3)}(t) & = - \int_D \int_0^{t-\zeta(x,z)-\sup\limits_{y\in D}\zeta(z,y)} \left(\frac{q(y)g(x,t-\tau;y)}{c_0^2(y)} - \frac{q(z)g(x,t-\tau;z)}{c_0^2(z)}\right)  u_{tt}(y,\tau) \, d\tau \, dy\\
  & \quad - \int_D \int_{t-\zeta(x,z)-\sup\limits_{y\in D}\zeta(z,y)}^{t-\zeta(x,y)} q(y)\,c_0^{-2}(y)g(x,t-\tau;y) u_{tt}(y,\tau)\, d\tau \, dy \\
  & \quad + \int_D \int_{t-\zeta(x,z)-\sup\limits_{y\in D}\zeta(z,y)}^{t-\zeta(x,z)} q(y)\,c_0^{-2}(y)g(x,t-\tau;z) u_{tt}(y,\tau)\, d\tau \, dy. 
\end{align*}
From the regularity properties of $g,\, q$ and $c_0$ alongwith the relation \eqref{v_H2}, we arrive at 
\begin{align}\label{reln_outside domain_2}
    \lvert {\textnormal{err}}_{(3)}(t) \rvert \preceq O(\epsilon^{-1}) \int_D \int_0^t |u_{tt}(y,\tau)| \, d\tau\, dy .
\end{align}
Using H\"older's inequality along with the relations \eqref{reln_outside domain_1}-\eqref{reln_outside domain_2} and the improved estimate \eqref{second bound}, we obtain
\begin{align}
    \label{final err_1} |{\textnormal{err}}_{(1)}(t)| & \preceq a^{-2}  \left( |D| \int_D |y-z|^2 \left|u_{tt}(y,t-\zeta(x,z))\right|^2 \, dy  \right)^{1/2} \preceq a^{2}, \\
    \nonumber |{\textnormal{err}}_{(3)}(t)| & \preceq a^{-2}  \left( |D| \int_D \int_0^t |y-z|^2 \left|u_{tt}(y,\tau)\right|^2 \, dy \, d\tau \right)^{1/2} \\
    & \preceq a^{-2}  \left( |D| \int_D |y-z|^2 \underset{0\le s\le t}{\textnormal{sup}}\, |\partial^2_t u(y,s)|^2\, dy \right)^{1/2} \preceq a^{2} \label{final err_2},
\end{align}
and in the same way
\begin{align}\label{final err_3}
    \nonumber |{\textnormal{err}}_{(2)}(t)| & \preceq a^{-2}  \left( |D| \int_D \left| \int_{\zeta(x,z)}^{\zeta(x,y)} \partial^3_t u(y,t-\kappa) \, d\kappa\right|^2 \, dy \right)^{1/2} \\
    \nonumber & \preceq a^{-2}  \left( |D| \int_D \left|\zeta(x,y)-\zeta(x,z)\right| \underset{0\le s\le t}{\textnormal{sup}}\, |\partial^3_t u(y,s)|^2\, dy \right)^{1/2} \\
    \nonumber & \preceq a^{-2}  \left( |D| \int_D |y-z| \underset{0\le s\le t}{\textnormal{sup}}\, |\partial^3_t u(y,s)|^2\, dy \right)^{1/2} \\
     & = a^2.
\end{align}
Putting together the error estimates \eqref{final err_1}-\eqref{final err_3} in \eqref{exapnsion_w}, we have
\begin{align}\label{outside domain_0}
    w(x,t) & = - \frac{\sigma(x,z)q(z)}{4 \pi c^3_0(z)\,|\zeta(x,\,z)|} \int_{D} u_{tt}(y,t-\zeta(x,z)) \, dy \, - \frac{q(z)}{c_0^2(z)} \int_{D} \int_{0}^{t} g(x,t-\tau;z)  u_{tt}(y,\tau)\, dy \, d\tau, \nonumber \\
    & \quad + O(a^2) ,
\end{align}
which holds pointwise for both space and time in $(\mathbb{R}^n \setminus \overline{D}) \times (0,T)$. Now the expansion \eqref{expansion in D} gives 
\begin{align}\label{inside domain_2nd derivative}
   u_{tt}(y,t) & = \sum_{n=1}^{\infty} \left(\frac{c_1}{\sqrt{\lambda_n}}\right)^3 e_n(y) \int_0^t \sin\left[\frac{c_1}{\sqrt{\lambda_n}}(t-\tau)\right]\, v(z,\tau)\, d\tau \int_D e_n(q) \, dq \nonumber \\
   & \quad + v(z,t) \, \sum_{n=1}^{\infty} \left(\frac{c_1}{\sqrt{\lambda_n}}\right)^2 e_n(y) \left(\int_D e_n(q) \, dq\right) \, + O(a^{2}),
\end{align}
to be understood point-wise in time and $L^2$-sense in $D$.
Inserting \eqref{inside domain_2nd derivative} into \eqref{outside domain_0}, we find 
\begin{align}\label{outside domain_1}
    w(x,t) & = - \sum_{n=1}^{\infty} \frac{\sigma(x,z)c_1}{4 \pi \lambda^{3/2}_n\, c_0(z) |\zeta(x,\,z)|} \left(\int_D e_n(y) \, dy\right)^2 \int_{0}^{t-\zeta(x,z)} \sin\left[\frac{c_1}{\sqrt{\lambda_n}}(t-\tau-\zeta(x,z))\right] \, v(z,\tau)\, d\tau \nonumber\\
    & \quad - \sum_{n=1}^{\infty} \frac{c_1}{\lambda_n^{3/2}} \left(\int_D e_n(y) \, dy\right)^2 \int_{0}^t  \int_0^{t -\tau}\sin\left(\frac{c_1}{\sqrt{\lambda_n}} (t -\tau-s) \right) g(x,s;z) v(z,\tau) \, ds\, d\tau \nonumber \\  
    & \quad - v(z,t-\zeta(x,z)) \, \sum_{n=1}^{\infty} \frac{\sigma(x,z)}{4 \pi \lambda_n\, c_0(z) |\zeta(x,\,z)|} \left(\int_D e_n(y) \, dy\right)^2 \nonumber \\
    & \quad - \sum_{n=1}^{\infty} \frac{1}{\lambda_n}  \left(\int_D e_n(y) \, dy\right)^2  \int_{0}^t g(x,t-\tau;z) v(z,\tau)\, d\tau \, + O(a^{2}),
\end{align}
which holds pointwise in both time and space variables in $( \overline{\Omega}\setminus V) \times (0,T)$. This completes the proof of Theorem \ref{outside domain_4}. Nonetheless, we continue the discussion to obtain a compact form of the expansion \eqref{outside domain_1} which will be used in Section \ref{recovery of speed and source}. Due to a change of variables and support condition of $v$ (in time), we observe that
\begin{align}\label{calc_1}
  \nonumber & \int_{0}^{t-\zeta(x,z)} \sin\left[\frac{c_1}{\sqrt{\lambda_n}}(t-\tau-\zeta(x,z))\right] \, v(z,\tau)\, d\tau \\
  \nonumber & \quad = \int_{\zeta(x,z)}^{t} \sin\left[\frac{c_1}{\sqrt{\lambda_n}}(t-\tau)\right] \, v(z,\tau-\zeta(x,z))\, d\tau \\
  & \quad = \int_{0}^{t} \sin\left[\frac{c_1}{\sqrt{\lambda_n}}(t-\tau)\right] \, v(z,\tau-\zeta(x,z))\, d\tau. 
\end{align}
Similarly, the support condition of $g$ (in time) implies
\begin{align}\label{calc_2}
   \nonumber & \int_{0}^t  \int_0^{t -\tau}\sin\left(\frac{c_1}{\sqrt{\lambda_n}} (t -\tau-s) \right) g(x,s;z) v(z,\tau) \, ds\, d\tau \\
   \nonumber & = \int_{0}^{t-\zeta(x,z)}  \int_0^{t -\tau}\sin\left(\frac{c_1}{\sqrt{\lambda_n}} (t -\tau-s) \right) g(x,s;z) v(z,\tau) \, ds\, d\tau \\
    \nonumber & = \int_{\zeta(x,z)}^{t} \int_0^{t -\tau + \zeta(x,z)}\sin\left(\frac{c_1}{\sqrt{\lambda_n}} (t -\tau+ \zeta(x,z)-s) \right) g(x,s;z) v(z,\tau-\zeta(x,z)) \, ds\, d\tau \\
    & = \int_{0}^{t} \int_0^{t -\tau + \zeta(x,z)}\sin\left(\frac{c_1}{\sqrt{\lambda_n}} (t -\tau+ \zeta(x,z)-s) \right) g(x,s;z) v(z,\tau-\zeta(x,z)) \, ds\, d\tau,
\end{align}
and also
\begin{align}\label{calc_3}
   \nonumber \int_{0}^t g(x,t-\tau;z) v(z,\tau)\, d\tau & = \int_{0}^{t-\zeta(x,z)} g(x,t-\tau;z) v(z,\tau)\, d\tau \\
   \nonumber & = \int_{\zeta(x,z)}^{t} g(x,t-\tau+\zeta(x,z);z) v(z,\tau-\zeta(x,z))\, d\tau \\
    & = \int_{0}^{t} g(x,t-\tau+\zeta(x,z);z) v(z,\tau-\zeta(x,z))\, d\tau.
\end{align}
For the ease of presentation, let us denote 
\begin{align}
   \label{func_1} a_n(x,t;z) & := - \frac{\sigma(x,z)c_1}{4 \pi \lambda^{3/2}_n\, c_0(z) |\zeta(x,\,z)|} \left(\int_D e_n(y) \, dy\right)^2 \sin\left( \frac{c_1}{\sqrt{\lambda_n}} t\right), \\
   \nonumber b_n(x,t;z) & := -  \frac{1}{\lambda_n} \left(\int_D e_n(y) \, dy\right)^2 \left[\frac{c_1}{\lambda_n^{1/2}}\int_0^{t+\zeta(x,z)} \sin\left(\frac{c_1}{\sqrt{\lambda_n}} (t+\zeta(x,z)-s) \right) g(x,s;z) \, ds \,\right. \\
   \label{func_2} & \hspace{8.5cm} \left. + \, g(x,t+\zeta(x,z);z) \right], \\
  \label{func_3}  \alpha(x,z) & := - \sum_{n=1}^{\infty} \frac{\sigma(x,z)}{4 \pi \lambda_n\, c_0(z) |\zeta(x,\,z)|} \left(\int_D e_n(y) \, dy\right)^2.
\end{align}
Combining the observations \eqref{calc_1}-\eqref{calc_3}, we can alternatively express \eqref{outside domain_1} as 
\begin{align}\label{outside domain_2}
   w(x,t) = \alpha(x,z) \, \tilde v(z,t) \, + \sum_{n=1}^{\infty} \int_{0}^{t} \left(a_n + b_n \right)(x,t-\tau;z) \, \tilde v(z,\tau) + O(a^{2}),
\end{align}
where $\tilde v(z,t) := v(z,t-\zeta(x,z))$.

%%%%%%%%%%%%%%%%%%%%%%%%%%%%%%%%%%%%%%%%%%%%%%%%%%%%%%%%%%%%%%%%%%%%%%%%%%%%%%%%%%%%%%%%
%%%%%%%%%%%%%%%%%%%%%%%%%%%%%%%%%%%%%%%%%%%%%%%%%%%%%%%%%%%%%%%%%%%%%%%%%%%%%%%%%%%%%%%%
%%%%%%%%%%%%%%%%%%%%%%%%%%%%%%%%%%%%%%%%%%%%%%%%%%%%%%%%%%%%%%%%%%%%%%%%%%%%%%%%%%%%%%%%
%%%%%%%%%%%%%%%%%%%%%%%%%%%%%%%%%%%%%%%%%%%%%%%%%%%%%%%%%%%%%%%%%%%%%%%%%%%%%%%%%%%%%%%%

\section{Proof of Theorem \ref{speed_source}}\label{recovery of speed and source}

We first show that the functions $a_n,\, b_n$ are summable (in appropriate sense) in time and $\alpha(x,z)$ is also finite (i.e. well-defined). For the reconstruction purpose and justification of the infinite series in \eqref{outside domain_2}, we need to know the eigenfunctions to the Newtonian operator $\mathcal{N}_D$. The  specific geometry of the droplet, which is spherical in our case, helps in this regard.  
\bigskip

Let us recall $D = z + a B(0,1)$ and the point that the asymptotic expansion of $w$ given in \eqref{outside domain_2} only involves those eigenfunctions of $\mathcal{N}_D$ which have non-zero averages. It is evident that, the translations of eigenfunctions become    
 eigenfunctions in that translated domain corresponding to the same eigenvalue. Therefore, if $\hat e_n$ denotes an eigenfunction to $\mathcal{N}_D$ corresponding to $\lambda_n$, then $\hat e_n$ is a suitable translation of $\tilde e_n$ which is related eigenfunction for the eigenvalue $\lambda_n$ to $\mathcal{N}_{B(0,a)}$. On that point, it suffices to discuss the eigensystem for $\mathcal{N}_{B(0,a)}$ denoted by $\{\lambda_n, \, \tilde e_n\}_{n\in\mathbb{N}}$. In this context, we follow the computations in \cite[Theorem 4.2]{K-S2011} to conclude 
\begin{align}\label{eigenfamily}
    \lambda^{-1}_{lj} = \left(\frac{\mu_j^{l+\frac{1}{2}}}{a} \right)^2, \, \textnormal{ and } \, \, \tilde e^m_{lj}(x) = \frac{1}{\sqrt{r}} J_{l+\frac{1}{2}}\left(\frac{r}{\sqrt{\lambda_{lj}}}\right) \, Y^m_{l}(\phi,\theta), \quad m\in\{-l,\cdots,l\}, 
\end{align}
for $(l,j) \in\mathbb{N}_0\times\mathbb{N}$, where $Y^m_{l}$ represents the spherical function denoted by 
\begin{align*}
    Y^m_l(\theta,\phi) = \begin{cases}
              P^m_l(\cos\theta)  \cos(m\phi), \quad m \in\{0,1,\cdots,l\}, \\
              P^{|m|}_l(\cos\theta) \sin(|m|\phi), \quad m \in\{-l,\cdot,-1\},
            \end{cases}
\end{align*}
and $\mu_j^{l+\frac{1}{2}}$ are roots to the equation
\begin{align*}
   (2l+1) J_{l+\frac{1}{2}}\left(\mu^{l+\frac{1}{2}}_j\right) + \frac{\mu^{l+\frac{1}{2}}_j}{2} \left[J_{l-\frac{1}{2}}\left(\mu^{l+\frac{1}{2}}_j\right) - J_{l+\frac{3}{2}}\left(\mu^{l+\frac{1}{2}}_j\right)\right] = 0.
\end{align*}
Here $P^m_l(\cdot)$ denotes the associated Legendre polynomial. For our specific interest, we need to consider those eigenfunctions having non-zero averages. In that context, we notice   
\begin{align*}
    \int_D \tilde e^m_{lj}(x) \, dx = 0 , \textnormal{ for } j\neq 0,
\end{align*}
since, $\int_{-1}^1 P^m_l(s) \, ds = 0$ for $m\in\{-l,\cdots,l\}$ and $l\neq 0$. As a consequence, we need to consider the eigenfamily \eqref{eigenfamily} only for $l=0$. Following the discussion in \cite{S-W2022}, we express this eigenfamily as
\begin{align}\label{eigensystem for ball}
   \lambda_n^{-1} = \frac{1}{a^2} \left(n\pi - \frac{\pi}{2} + \gamma_n\right)^2, \quad
   \tilde e_n (x)= \frac{1}{\sqrt{|x|}} J_{\frac{1}{2}}\left(\frac{|x|}{\sqrt{\lambda_n}}\right) = \sqrt{\frac{2\sqrt\lambda_n}{\pi}}\,\frac{1}{|x|}\, \sin\left(\frac{|x|}{\sqrt{\lambda_n}}\right), \, n\in \mathbb{N} ,  
\end{align}
where $\{\gamma_n\}_{n\in\mathbb{N}}$ is a monotonically decreasing sequence converging to zero.

For these special eigenfamily, we perform explicit computations in order to justify the existence of $\alpha(x,z)$ and derive asymptotic behaviors of the functions for $a_n$ and $b_n$. To this aim, let us set 
\begin{equation*}
m_n := \lambda_n^{-1/2}a = n\pi - \frac{\pi}{2} + \gamma_n,
\end{equation*}
and hence, we have $\lim\limits_{n\to\infty} m_n = \infty$. An important relation  in this context is 
\begin{align}\label{tan reln}
   \tan m_n = - 2m_n,
\end{align}
which has been already derived in \cite{S-W2022}. Therefore, we deduce 
\begin{align}
    \nonumber \int_D \hat e_n(y) \,dy  = \int_{B(0,a)} \tilde e_n(x)\, dx & = \sqrt{8\pi^3\,\sqrt\lambda_n} \, \int_0^a r \sin\left(\frac{r}{\sqrt \lambda_n}\right) \,dr \\
    & = \sqrt{8\pi^3} \lambda_n^{\frac{5}{4}}\, \int_0^{m_n} r \sin r \,dr \nonumber \\
    \nonumber & = \sqrt{8\pi^3} \lambda_n^{\frac{5}{4}} \, \left(\sin m_n - m_n \cos m_n\right) \\
   \nonumber & = -3 \sqrt{2\pi^3} \lambda_n^{\frac{5}{4}} \, \sin m_n \\
   \label{integral of e_n} & = (-1)^n 3 \sqrt{2\pi^3} \, \lambda_n^{\frac{5}{4}} \, \cos \gamma_n,
  \end{align}
where we utilized \eqref{tan reln} in last line. Similarly, we make use of \eqref{tan reln} and compute 
\begin{align}
  \nonumber \|\hat e_n\|^2_{L^2(D)} = 4\pi \sqrt{\lambda_n} \int_0^a \sin^2\left( \frac{r}{\sqrt\lambda_n}\right)dr & = 2\pi\lambda_n \left( m_n - \sin m_n \cos m_n\right) \\
  \label{a part of integral of e_n} & = 2\pi \lambda_n \left(m_n + \frac{\cos^2\gamma_n}{m_n}\right).
\end{align}
As we have considered orthonormal family in the expansion \eqref{exapnsion_w}, we normalize $\hat e_n$ and set 
$$ e_n(x) = \frac{1}{\|\hat e_n\|_{L^2}} \, \hat e_n(x), \quad \textnormal{ for } x \in D, \, n\in\mathbb{N}. $$
Combining \eqref{integral of e_n} and \eqref{a part of integral of e_n} with the fact that, $\cos\gamma_n \to 1$, as $n\to\infty$, we  notice that 
\begin{align}\label{int of e_n}
   \left\vert \int_D e_n(x) \, dx\right\vert \simeq \lambda_n \simeq n^{-2}, \textnormal{ for } n>>1,  
\end{align}
implying convergence of the series
\begin{align*}
    \sum_{n=1}^\infty \frac{1}{\lambda_n}\,\left(\int_D e_n(x) \, dx\right)^2 .
\end{align*}
It proves that $\alpha(x,z)$ and the fourth infinite series in the right hand side of \eqref{outside domain_1} make sense. Now, to justify the convergence of first two infinite series of \eqref{outside domain_1}, we argue as follows. Without loss of generality, we may assume that $ c_1\,T < a$. This is due to the fact that $c_1$ being related to the injected droplet can be adjusted accordingly. In light of this, Theorem 1.2 of \cite{S-W2022} yields that $\left\{\sin\left(\frac{c_1\, t}{\sqrt{\lambda_n}}\right)\right\}_{n\in\mathbb{N}}$ defines a Riesz basis in $L^2\left(-\frac{a}{c_1}, \, \frac{a}{c_1}\right)$ and hence in $L^2(0,T)$ too. Due to this and the asymptotic relation \eqref{int of e_n},  we notice that the function
\begin{align*}
    S(t)=\sum_{n=1}^\infty\, a_n(x,t;z), \quad \textnormal{ for a.e. } t \in (0,T),
\end{align*}
defines an $L^2$ function. Therefore, the first infinite series in \eqref{outside domain_1} is nothing but the convolution (in time) between $v(z, \cdot)$ and $S(\cdot)$. A similar observation holds for the second infinite series in \eqref{outside domain_1}. Summarizing all these arguments, we are able to give sense to the dominant term of the wave-field $w(x,t)$ through an elementary calculations which proves to be useful next.

\bigskip

Now, we are interested in recovering the initial wave-field $v$ from the dominating part of the measurement. In other words, we recast \eqref{outside domain_2} as
\begin{align}\label{int form_data}
   w(x,t) = \mathcal{A}\tilde v(z,t) + O(a^{2})
\end{align}
and wish to invert the operator $\mathcal{A}:= \alpha I + \mathcal{K}$ in $L^2(0,T)$ where
\begin{align}\label{defn_conv op}
\begin{cases}
   \mathcal{K}f(t) =  \int_0^t K(t-s) f(s) \, ds, \\
   \alpha = \alpha(x,z), \quad  K(t) = \sum_{n=1}^{\infty} (a_n + b_n)(x,t;z).
\end{cases}
\end{align}
 The kernel of the convolution operator $\mathcal{K}$ being an $L^2$ function implies that $\mathcal{K}$ is a Hilbert-Schmidt operator and hence compact in $L^2(0,T)$. Now we apply Fredholm alternative to infer that $\mathcal{A}$ is invertible in  $L^2(0,T)$ iff $\mathcal{K}$ does not admit $-\alpha$ as an eigenvalue. In  the latter case, there exists $f\in L^2(0,T)$ such that $f = -\frac{1}{\alpha}\mathcal{K}f$, implying $f\in L^\infty(0,T)$, since $\mathcal{K}$ has $L^2$ kernel. Therefore, it suffices to establish invertibility of $\mathcal{A}$ in $L^\infty(0,T)$.

To show this, we first see that the definition of $\mathcal{K}$ in \eqref{defn_conv op} yields that 
\begin{align}\label{L^2 induction_0}
    |\mathcal{K}^nf(t)| \le \frac{\|K\|^n \, \|f\|_{L^\infty}}{\sqrt{n!}}\, t^{\frac{n}{2}}, \quad \textnormal{ for a.e. } t \in (0,T), \, n \in \mathbb{N}   , 
\end{align}
and hence, it is clear that
\begin{align}\label{L^2 induction}
     \|\mathcal{K}^n\|_{L^\infty(0,T)} \le \frac{(\|K\| \sqrt{T})^n}{\sqrt{n!}}\,  , \quad n \in \mathbb{N}.
\end{align}
Using H\"older's inequality, we verify \eqref{L^2 induction_0} by an induction argument. Showing \eqref{L^2 induction_0} for $n=1$ and a.e. $t\in(0,T)$ is not hard, since 
\begin{align*}
   \left|\mathcal{K}f(t)\right| = \left| \int_0^t K(t-s)\, f(s)\, ds\right| \le \|f\|_\infty \left( \int_0^t |K(s)|^2 \, ds\right)^{1/2} \, \left( \int_0^t \,ds\right)^{1/2} \le \|f\|_\infty \|K\| \sqrt{t}.
\end{align*}
Let us assume that \eqref{L^2 induction_0} holds for $n=m$ and proceed to prove it for $n=m+1$. For a.e. $t\in(0,T)$, 
\begin{align*}
    \left|\mathcal{K}^{m+1} f(t)\right| =  \left|\int_0^t K(t-s)\,\mathcal{K}^{m} f(s) \, ds \right| & \le \left( \int_0^t |K(s)|^2 \, ds\right)^{1/2} \, \left( \int_0^t |\mathcal{K}^{m} f(s)|^2\,ds\right)^{1/2} \\
    & \le \frac{\|K\|^{n+1}\, \|f\|_\infty}{\sqrt{n!}} \,\left( \int_0^t s^n\,ds\right)^{1/2} \\
    & = \frac{\|K\|^{n+1}\, \|f\|_\infty}{\sqrt{(n+1)!}}\, \,t^{\frac{n+1}{2}}.
\end{align*}
Therefore, the inequality \eqref{L^2 induction_0} is also true for $n=m+1$ proving our claim. 
As a consequence of \eqref{L^2 induction} and Gelfand's theorem on spectral radius, we deduce
\begin{align}\label{Gelfand}
  \rho(\mathcal{K}) = \liminf\limits_{n\to\infty} \|\mathcal{K}^n\|_{L^\infty(0,T)}^{1/n} \le  \|K\| \sqrt{T}\, \lim\limits_{n\to\infty} \,\frac{1}{\sqrt[n]{n!}}
\end{align}
where $\rho(\mathcal{K})$ denotes the spectral radius of $\mathcal{K}$ in $L^\infty$ topology. 
Define $x_n := \frac{1}{\sqrt{n!}}$ and we see that $\lim\limits_{n\to\infty} x_n^{1/n} = 0$, since 
$
  \frac{x_{n+1}}{x_n} = \frac{1}{\sqrt{n+1}} \longrightarrow 0,  \textrm{ as } n \to \infty.
$
From \eqref{Gelfand}, we conclude that $\rho(\mathcal{K})=0$. Hence we can invert the operator $\mathcal{A}$ in $L^\infty(0,T)$ and then in $L^2(0,T)$ as mentioned before.  Let us introduce the inverse of $\mathcal{A}$ in $L^2(0,T)$ as
\begin{align}\label{inverse_series expansion}
    \mathbb{A} := (\alpha I + \mathcal{K})^{-1} 
     = \sum_{n=0}^{\infty} \alpha^{-n-1}  \mathcal{K}^{n}.
\end{align} 
 From the discussion in Section \ref{Preliminaries}, it is clear that the term $\alpha$ and the function $K(\cdot)$ are multiplied with the scaling parameter $a$, implying that $\mathbb{A}w(x,\cdot) = O(1)$ to be understood pointwise in time variable. As a consequence, we express \eqref{int form_data} as 
\begin{align*}
    v(z,t) = \mathbb{A}w(x,\cdot)(t+\zeta(x,z)) + O(a),
\end{align*}
which proves the first part of Theorem \ref{speed_source}.

\bigskip

Now we want to recover the function $\zeta(x,z)$. Let us first highlight an important feature of the dominating term of the asymptotic expansion \eqref{outside domain_2} i.e. $\mathcal{A}\tilde v(z,\cdot)$, which indicates that 
\begin{align*}
\mathcal{A}\tilde v(z,t) = 0, \quad \textrm{ for } t<\zeta(x,z).
\end{align*}
It is due to the support condition of $v(z,\cdot)$. Furthermore, the graph of $w(x,\cdot)$ vanishes trivially before the time level $t=\inf\limits_{y\in D}\zeta(x,y)= \zeta(x,z) + O(\epsilon)$. To determine $\zeta(x,z)$, it suffices to show that $w(x,\cdot)$ is not identically zero after $t=\zeta(x,z)$. However, we are only able to prove this for $\mathcal{A}\tilde v(z,\cdot)$ which is the dominant term of $w(x,\cdot)$. Hence, we can recover $\zeta(x,z)$ upto an error. From the invertibility of $\mathcal{A}$, the  question whether $\mathcal{A}\tilde v(z,\cdot)$ vanishes identically after $t=\zeta(x,z)$ can be rephrased by saying whether for every $\delta>0$, we have
\begin{align*}
   \tilde v(z,t^*_{\delta,z}) \neq 0, \quad \textrm{ for some } t^*_{\delta,z} \in \left(\zeta(x,z) ,\, \zeta(x,z)+\delta \right),
\end{align*}
or equivalently,
\begin{align}\label{nonzero_statement_1}
   v(z,t_{\delta,z}) \neq 0, \quad \textrm{ for some } t_{\delta,z} \in (0,\, \delta).
\end{align}
We can not always show \eqref{nonzero_statement_1} since 
\begin{align}\label{rep formula_v}
  \nonumber  v(z,t) & = \int_\mathbb{R} \int_{\mathbb{R}^n} G(z,t;y, \tau) J(y,\tau) \, dy d\tau \\
   & = \frac{1}{4\pi} \int_\Omega \frac{\sigma(z,y)}{c_0(y)\,|\zeta(z,\,y)|} J(y,t-\zeta(z,y)) \, dy \, + \int_{\Omega} \int_{0}^t g(z,t-\tau;y) J(y,\tau) \, dy d\tau,
\end{align}
and thus the non-vanishing of source function $J$ (in time) determines whether the wave-field $v$ should be zero or not.  
Hence the assumption that $J(x,\cdot)$ does not trivially vanish for a.e. $x\in\Omega$ is quite reasonable and this is what has been assumed in Theorem \ref{speed_source}. Because of the regularity of the wave speed $c_0(\cdot)$ and $\zeta(x,\cdot)$, we remark that finding the function $\zeta(x,z)$ for $z$ varying in some dense set of $\Omega$ is enough to determine $c_0$ in the whole of $\Omega$. Keeping this in mind, we consider 
\begin{align}\label{nonzero_statement_2}
   \mathcal{W}:= \{ z\in \Omega; \, \forall \delta>0, \, \exists \, t_{\delta,z} \in (0,\, \delta) \textrm{ such that } v(z,t_{\delta,z}) \neq 0\},
\end{align}
and claim that $\mathcal{W}$ is dense in $\Omega$. If not, then we have a closed ball, say $B_* \subseteq \Omega$, such that $\forall b \in B_*$, there exists $t_b >0$ and 
\begin{align*}
    v(b,t)= 0, \, \textrm{ for } t \in (0, \, t_b).
\end{align*}
Using the compactness of $B_*$, we take $t_* = \inf\limits_{b\in B_*} t_b > 0$ and then we have 
\begin{align*}
    v(b,t) = 0 ,\quad \textrm{ for all } b\in B_* \textrm{ and } t < t_*,
\end{align*}
which implies $\left. J(\cdot,\cdot)\right|_{B_* \times (0,t_*)} \equiv 0$ contradicting our assumption on $J$. By virtue of the space regularity of $v$ as $v(\cdot,t) \in C^{0,\frac{1}{2}}(\Omega)$, we also find that $\mathcal{W}$ is an open set of $\Omega$. Therefore we conclude that $\mathcal{W}$ is an open dense set of $\Omega$. We summarize the above discussion by saying that the graph of $\mathcal{A}\tilde{v}(z,\cdot)$ experiences jump at $t=\zeta(x,z)$ under the assumption that the source term $J$ is non-vanishing near $t=0$, where $z$ is arbitrarily chosen from the dense set $\mathcal{W}$. Hence, we prove that the function $\tilde v(z,\cdot)$ and therefore $\mathcal{A}\tilde v(z,\cdot)$ do not identically vanish after $t=\zeta(x,z)$. 
\bigskip

Now we study the propagation of error while determining the function $\zeta(x,\cdot)$ from the measurement $w(x,\cdot)$. We see from \eqref{outside domain_2} that all the coefficients in the linear operator $\mathcal{A}$ depend continuously w.r.t the variable $z \in \Omega$. To underline this dependence, we denote the linear operator $\mathcal{A}$ by $\mathcal{A}_z$ and then we notice that the map $ p \to \mathcal{A}_p$ is continuous with 
\begin{align}\label{Holder cont_op}
  \|\mathcal{A}_p - \mathcal{A}_q\|_{L^2(0,t)} \preceq |p-q|, \quad \forall t \in (0,T).
\end{align}
In pursuit of recovering the wave speed $c_0(z)$ for some fixed $z \in \Omega$, we utilize the denseness of $\mathcal{W}$ to ensure that there exists $z_\epsilon \in \mathcal{W}$ such that $|z_\epsilon - z| < \epsilon$. We observe
\begin{align}\label{speed deter_approax}
   w(x,t) & = \mathcal{A}_z \tilde v(z,t) + O(a^{5/2}) \nonumber \\
    & = \mathcal{A}_{z_\epsilon}\tilde v(z_\epsilon,t) + \left(\mathcal{A}_z - \mathcal{A}_{z_\epsilon}\right)\tilde v(z,t) + \mathcal{A}_{z_\epsilon} \left(v(z_\epsilon,t) - v(z,t)\right) + O(a^{2}),
\end{align}
and due to Lemma \ref{a-priori estimate} applied to $v$, we have $v(\cdot,t) \in C^{0,\frac{1}{2}}(\Omega)$ for all $t\in [0,T]$. Thus, we have
\begin{align}\label{Holder cont}
   |\tilde v(z_\epsilon,t) - \tilde v(z,t)| \preceq |z_\epsilon - z|^{\frac{1}{2}} \preceq \sqrt{\epsilon}.
\end{align}
 We recall that the definition of $\mathcal{A}$ involves the term $\alpha$ and the function $K(\cdot)$ which are multiplied by the scaling parameter $a$. Inserting the observations \eqref{Holder cont_op} and \eqref{Holder cont} into \eqref{speed deter_approax}, we  write
\begin{align}\label{error_w}
     w(x,t) = \mathcal{A}_{z_{\epsilon}} \tilde v(z_{\epsilon},t) + O(a\epsilon)+ O(a\sqrt{\epsilon}) + O(a^{2}),
\end{align}
which holds point-wise in both space and time variables. Now taking $\epsilon = a^2$ in \eqref{error_w}, we rewrite
\begin{align*}
    w(x,t) = \mathcal{A}_{z_{\epsilon}} \tilde v(z_{\epsilon},t) + O(a^{2}).
\end{align*}
In conclusion, we recover the term $\mathcal{A}_{z_{\epsilon}}\tilde v(z_{\epsilon},\cdot)$ for $t\in (0,T)$ from the measurement $w(x,t)\vert_{t\in (0,T)}$ modulo an error term of order $a^{2}$ in point-wise sense w.r.t time variable. Following our earlier discussion, we can identify the time level where the jump happens in the graph of $\mathcal{A}_{z_{\epsilon}}\tilde v(z_{\epsilon},\cdot)$, which is the time level $t=\zeta(x,z_{\epsilon})$. Using the regularity of $\zeta(x,\cdot)$, we can say that $\zeta(x,z)$ can be recovered from the measurement $w(x,\cdot)$ with an error of the order $O(a^2)$. Injecting the droplet at different points $z \in \Omega$, the functions $\zeta(x,z)$ and $v(z,t)$ are recovered. Then we perform numerical differentiation in \eqref{Eikonal} and \eqref{model_v} to first recover the wave speed $c_0(\cdot)$ and the source function $J(\cdot,\cdot)$ respectively. 

\bigskip

\section{Discussions on Theorem \ref{speed_source}} \label{remarks}

 We have shown that we can recover simultaneously the wave speed and source function. Here, we highlight two related aspects we consider to be worth discussing. 
 \begin{enumerate}
 \item The first one is that we can not determine the function $\zeta(x,z)$ for all $z\in \Omega$, directly from the graph of $v(z,\cdot)$. The possibility that the function $v(z,\cdot)$ is trivially zero for some $z\in \Omega$ can not be avoided even when the source function $J(z,\cdot)$ is non-vanishing for all $z\in \Omega$. Let us provide a small example justifying the last claim for the simple case when $c_0 \equiv 1$. For $r_1, r_2 \in \mathbb{R}$ satisfying $r_1 > r_2 >0$, consider $\Omega:= B(0,r_1) \subset \mathbb{R}^3$ and $u(x,t) = (|x|^2 - r^2_2)t^2$. Then $u$ satisfies
 \begin{align*}
     \begin{cases}
         (\partial_t^2 - \Delta) u(x,t) = J(x,t):=  2(|x|^2 - r^2_2 - 3t^2),    \\
         u(x,0) = u_t(x,0) = 0.
     \end{cases}
 \end{align*}
 We notice that the function $J(x,\cdot)$ is not identically zero for all $x$, although $u(x,\cdot)\rvert_{x \in \partial B(0,r_2)} \equiv 0$. Therefore we expect the non-vanishing of $v(z,\cdot)$ to happen at most in a dense set of $\Omega$, which has been shown in Section \ref{recovery of speed and source}. 
 
  \item The second one is related to the sharpness of the jump. Indeed, it is naturally desirable to detect a sharp jump in the graph of $\mathcal{A}v(z,\cdot)$ at $t=0$ rather than just having $v(z,\cdot)$ as a non-vanishing function for $t>0$. We discuss this aspect here. From this moment on, we assume that the source function $J$ is smooth and has some non-zero time-derivative on the plane $t=0$. To acknowledge the condition of $J$ as mentioned in Lemma \ref{a-priori estimate}, we need to assume 
\begin{align}\label{assumption on J_1}
\partial^{p+1}_t J(y,0) \neq 0, \textrm{ and } \partial^k_t J(y,0) = 0, \ \forall k \in \{0,1,..,p\}, \ \forall y \in \Omega.  
\end{align}
For our arguments in Section \ref{sec_asy}, we took $p=13$ and continue with this choice. Without loss of generality, let us assume that, $\partial^{p+1}_t J(z,0) > 0$. Due to smoothness of $J$, we can find some $ \alpha,\ r_0 >0$ and $\delta>0$ so that 
\begin{align}\label{assumption on J_2}
    \partial^{p+1}_t J(y,t) \ge \alpha > 0, \textrm{ for } (y,t) \in \mathcal{B}_{z,\delta}:= B(z,r_0) \times [0,\delta]. 
\end{align}
 We use the relations \eqref{assumption on J_1} and \eqref{assumption on J_2} in  Taylor's expansion for $J(y, \cdot)$ at $t=0$
\begin{align}\label{Taylor}
   J(y,t) = \sum_{k=1}^{p} \frac{\partial^{k}_t J(y,0)}{k!} t^k + \frac{\partial^{p+1}_t J(y,\theta t)}{(p+1)!} t^{p+1} = \frac{ \partial^{p}_t J(y,0)}{p!} t^p + \frac{\partial^{p+1}_t J(y,\theta t)}{(p+1)!} t^{p+1}.
\end{align}
Let us further choose $\delta$ sufficiently small so that \eqref{Taylor} yields for some constants $c_1 , c_2 >0$ 
\begin{align}\label{ineq for J}
    c_1 t^p  \le J(y,t) \le c_2 t^p ,\quad \forall (y,t) \in \mathcal{B}_{z,\delta}.
\end{align}
To detect the sharpness in the graph of $v(z, \cdot)$ near $t=0$, it suffices to restrict our discussion only for small $t$. For this reason, let us fix $z \in \Omega$ and $t>0$ to define the region 
\begin{align*}
    \Omega_t : = \{y \in \Omega; \ \zeta(z,y) < t\}.
\end{align*}
We have already seen that $\zeta(z,y)$ denotes the Reimannian distance between $z$ and $y$ w.r.t the metric \eqref{metric}. From \eqref{v_H2}, we notice that there exists constants $d_1, d_2 >0$ satisfying 
\begin{align*}
    d_1 |y-z| < \zeta(z,y) < d_2 |y-z|, \ \textrm{ when } y \textrm{ is near } z,
\end{align*}
which in turn implies 
\begin{align}\label{metric comparison}
    B(z,t/{d_2}) \subseteq \Omega_t \subseteq B(z,t/{d_1}).
\end{align}
From the causality of $J$ and support condition on $g$, we find
\begin{align*}
    v(z,t) = \frac{1}{4\pi} \int_{\Omega_t} \frac{\sigma(z,y)}{c_0(y)\,|\zeta(z,\,y)|} J(y,t-\zeta(z,y)) \, dy \, + \int_{0}^t \int_{\Omega_{t-\tau}} g(z,t-\tau;y) J(y,\tau) \, dy d\tau.
\end{align*}
For $t< \delta$, we use continuity of $\sigma(z,\cdot)$ and $c^{-1}_0(\cdot)$ along with the relations \eqref{ineq for J} and \eqref{metric comparison} to obtain 
\begin{align}\label{1st term}
   \nonumber   \frac{1}{4\pi} \int_{\Omega_t} \frac{\sigma(z,y)}{c_0(y)\,|\zeta(z,\,y)|} J(y,t-\zeta(z,y)) \, dy & \succeq \int_{B(z,t/{d_2})} \frac{\left(t- \zeta(z,y)\right)^p}{\zeta(z,y)} \, dy \\
    \nonumber  & \succeq \int_{B(z,t/{d_2})} \frac{\left(t- d_2 |y-z|\right)^p}{|y-z|} \, dy \\
  & \succeq \int_{B(z,t)} \frac{\left(t- |y-z|\right)^p}{|y-z|} \, dy \nonumber \\
  & \succeq \int_0^t r(t-r)^p \, dr \simeq t^{p+2}
\end{align}
and 
\begin{align}\label{2nd term}
        \left|\int_{0}^t \int_{\Omega_{t-\tau}} g(z,t-\tau;y) J(y,\tau) \, dy d\tau \right| \preceq \int_{0}^{t} \int_{B(z,\frac{t-\tau}{d_1})} \tau^{p} \, d y\, d\tau \preceq \int_{0}^{t} (t-\tau)^3 \tau^{p} d\tau \preceq t^{p+4}.
\end{align}
From \eqref{1st term} and \eqref{2nd term}, we obtain $v(z,t) \succeq t^{p+2}$ for $t \in (0,\delta)$. In the same way, we also have the reverse inequality $v(z,t) \preceq t^{p+2}$ for $t \in (0,\delta)$, implying 
\begin{align*}
    v(z,t) \simeq t^{p+2}, \quad \textnormal{ for } t \in (0,\delta).
\end{align*}
For $\delta>0$ small, we find
\begin{align*}
    \mathcal{A}v(z,t) = \alpha v(z,t) + \int_0^t K(t-\tau) v(z,\tau) \, d\tau \ge \alpha t^{p+2} - aC\, \int_0^t \tau^{p+2} \, d\tau \succeq at^{p+2}, \quad t\in (0,\delta).
\end{align*}
Taking this into consideration, we notice 
\begin{align*}
w(x,t) = \mathcal{A}\tilde{v}(z,t) + O(a^2) \succeq a(t-\zeta(x,z))^{p+2} - Ca^2, \quad \textnormal{ for } t\in \left( \zeta(x,z), \zeta(x,z) +\delta \right),
\end{align*}
and as a result, one always experiences jump in the graph of $w(x,\cdot)$ at  time level
\begin{equation}\label{sharpness}
t_* = \zeta(x,z) + a^{\frac{1-\epsilon}{p+2}} \mbox{ where } \epsilon \in (0,1).
\end{equation}
Therefore, we can detect the travel-time level at an error of the order $a^{\frac{1-\epsilon}{p+2}}$, $a \ll 1$.
\end{enumerate}

%%%%%%%%%%%%%%%%%%%%%%%%%%%%%%%%%%%%%%%%%%%%%%%%%%%%%%%%%%%%%%%%%%%%%%%%%%%%%%%%%%%%%%%%
%%%%%%%%%%%%%%%%%%%%%%%%%%%%%%%%%%%%%%%%%%%%%%%%%%%%%%%%%%%%%%%%%%%%%%%%%%%%%%%%%%%%%%%%
%%%%%%%%%%%%%%%%%%%%%%%%%%%%%%%%%%%%%%%%%%%%%%%%%%%%%%%%%%%%%%%%%%%%%%%%%%%%%%%%%%%%%%%%
%%%%%%%%%%%%%%%%%%%%%%%%%%%%%%%%%%%%%%%%%%%%%%%%%%%%%%%%%%%%%%%%%%%%%%%%%%%%%%%%%%%%%%%%

\bigskip
{\bf Acknowledgement:} This work is supported by National Natural Science Foundation of China (Nos. 12071072, 11971104) and the Austrian Science Fund (FWF): P 30756-NBL.

\end{document}